\documentclass[11 pt]{article}
\usepackage{style}
\title{\LARGE\bfseries Maximizing the Value of Predictions in Control:\\
Accuracy Is Not Enough}
\author{
Yiheng Lin\textsuperscript{1},
Christopher Yeh\textsuperscript{1},
Zaiwei Chen\textsuperscript{2}, and
Adam Wierman\textsuperscript{1}\\
{\small\textsuperscript{1}\textit{Caltech, Department of Computing and Mathematical Sciences, Pasadena, CA, USA}}\\
{\small\textit{\href{mailto:yihengl@caltech.edu}{yihengl@caltech.edu}, \href{mailto:cyeh@caltech.edu}{cyeh@caltech.edu}, \href{mailto:adamw@caltech.edu}{adamw@caltech.edu}}}\\
{\small\textsuperscript{2}\textit{Purdue, School of Industrial Engineering, West Lafayette, IN, USA,} \href{mailto:chen5252@purdue.edu}{\textit{chen5252@purdue.edu}}}
}

\date{\vspace{-0.4 in}}
\begin{document}
\maketitle

\begin{abstract}
We study the value of stochastic predictions in online optimal control with random disturbances. Prior work provides performance guarantees based on prediction error but ignores the stochastic dependence between predictions and disturbances. We introduce a general framework modeling their joint distribution and define ``prediction power'' as the control cost improvement from the optimal use of predictions compared to ignoring the predictions. In the time-varying Linear Quadratic Regulator (LQR) setting, we derive a closed-form expression for prediction power and discuss its mismatch with prediction accuracy and connection with online policy optimization. To extend beyond LQR, we study general dynamics and costs. We establish a lower bound of prediction power under two sufficient conditions that generalize the properties of the LQR setting, characterizing the fundamental benefit of incorporating stochastic predictions. We apply this lower bound to non-quadratic costs and show that even weakly dependent predictions yield significant performance gains.
\end{abstract}

\section{Introduction}\label{sec:intro}
Understanding the benefits of predictions in control has received significant attention recently \cite{yu2020power,yu2022competitive,lin2021perturbation,lin2022bounded,zhang2021regret}. In this work, we study a class of discrete-time online optimal control problems in general time-varying systems, where random disturbances $W_t$ affect state transitions. The agent leverages a prediction vector containing information about future disturbances to minimize the expected total cost over a finite horizon $T$. To study the impact of using predictions, a fundamental question is how to model disturbances and their relationship to predictions. Prior works adopt different modeling approaches \cite{yu2020power,zhang2021regret,chen2015online}, each with distinct strengths and limitations.

A common paradigm assumes perfect predictions over a finite horizon $k$, yielding an elegant characterization of ``prediction power'' that improves with larger $k$. Under this model, predictions exactly reveal future disturbances $W_t, \ldots, W_{t+k-1}$. \cite{yu2020power} shows how prediction power grows with $k$ in the LQR setting, and subsequent work extends this result to time-varying systems \cite{lin2021perturbation}. As a result, the marginal benefit of one additional prediction decays exponentially with $k$, offering insight into how to select $k$. However, longer-horizon predictions are more costly and less accurate, and real-world predictions are rarely perfect \cite{chen2024soda}, making this idealized setting challenging in practice.

A natural extension of accurate predictions is to consider bounded prediction errors, which better captures practical challenges. Specifically, prediction errors measure the distance between the predicted and actual disturbances, and the resulting cost bounds depend on these errors \cite{yu2022competitive,zhang2021regret,lin2022bounded}. This extension recovers the perfect predictions setting when errors shrink to zero. However, it can be overly pessimistic because, for any predictor, the same performance bound must also apply to an adversary that generates the worst prediction sequence to penalize the predictive policy subject to the same error bound. It overlooks stochastic dependencies between predictions and disturbances that can be valuable for improving control costs.

In this work, we propose a general stochastic model that captures the distributional dependencies between predictions and disturbances, without restricting prediction targets, horizon length, or requiring strict error bounds. Compared with previous stochastic methods \cite{chen2015online,chen2016using}, our approach further relaxes problem-specific assumptions and directly focuses on the incremental benefit of predictions. Such benefits can be subtle—often overlooked by classical metrics like regret or competitive ratio. To capture them, we define \emph{prediction power} as the improvement in expected total cost when predictions are fully exploited, which builds on and generalizes the notion from \cite{yu2020power}. Our framework thus characterizes when and why predictions significantly boost online control performance.

\textbf{Contributions.} We introduce a general stochastic model (\Cref{def:predictions}) that describes how disturbances relate to all candidate predictors. We then define \emph{prediction power} (\Cref{def:prediction-power}) which quantifies the incremental control-cost improvement gained by fully leveraging these predictions. To illustrate this concept, we derive an exact expression for prediction power in the benchmark setting of time-varying linear quadratic regulator (LQR) control (\Cref{thm:pred-power-LTV}). Using this closed-form formula, we provide examples (\textit{e.g.}, \Cref{example:one-step-dim-rotation}) that illustrate why analyzing prediction accuracy is insufficient---improving prediction accuracy may not always improve prediction power. Finally, we demonstrate the connection between prediction power and online policy optimization (\Cref{example:pred-power-and-policy-opt}), highlighting how practical algorithms can attain (a portion of) the maximum potential.

We extend our analysis of prediction power beyond the LQR setting. This generalization poses significant challenges due to the lack of closed-form expressions for the optimal policy. Building on insights from the LQR analysis, we identify two key structural conditions: a quadratic growth condition on the optimal Q-function (\Cref{cond:Q-function-quadratic-growth}) and a positive semi-definite covariance condition on the optimal policy's actions (\Cref{cond:optimal-policy-psd-covariance}). These conditions are sufficient to derive a general lower bound on prediction power, formalized in \Cref{thm:Bellman-improvement-not-worse}. We apply this result to the setting of time-varying linear dynamics with non-quadratic cost functions. Under assumptions on costs and on the joint distribution of predictions and disturbances, we establish a lower bound on prediction power (\Cref{coro:one-step-pred-power-lower}), demonstrating that even weak predictions can yield strict performance gains.

\textbf{Related Literature.} Our work is closely related to the line of works on using predictions in online control. Our prediction power is inspired by \cite{yu2020power}, which defines the prediction power as the maximum control cost improvement enabled by $k$ steps of accurate predictions in the time-invariant LQR setting. Compared with \cite{yu2020power}, we extend the notion of prediction power to allow general dependencies between predictions and disturbances, and we consider more general dynamics/costs (\Cref{sec:main_results}). Rather than focusing on the prediction power, many works study the power of a certain policy class such as MPC \citep{yu2022competitive,lin2021perturbation,zhang2021regret,lin2022bounded}, Averaging Fixed Horizon Control \citep{chen2015online,chen2016using}, Receding Horizon Gradient Descent \citep{li2018using,li2019online}, and others \citep{lin2020online}.
While one can say the power of (generalized) MPC equals to the prediction power in the LQR setting \citep{yu2020power} (Section \ref{sec:LQR}), they are not the same in general (see Appendix \ref{appendix:MPC-counterexample}).

Our work is, in part, motivated by both empirical and theoretical findings in the decision-focused learning (DFL) literature that prediction models with the same prediction accuracy may have very different control costs (see \cite{mandi2024decision} for a recent survey). Research on DFL typically considers predictions given as point estimates of some uncertain input to decision-makers modeled as optimization problems, such as stochastic optimization (\cite{donti2017task}), linear programs (\cite{elmachtoub2022smart}), or model predictive control (\cite{amos2018differentiable}), although more recent works have started exploring other forms of predictions such as prediction sets (\cite{yeh2024end,wang2024learning}). In contrast, our work does not require any particular form of decision-maker; instead, our main result characterizes the benefit of optimally leveraging predictions, for whatever form an optimal controller may take. Whereas DFL aims to design procedures for training prediction models that reduce downstream control costs, our work studies a more fundamental question about how much performance gain is achievable with better predictions.

\section{Problem Setting}\label{sec:setting}
We consider a finite-horizon discrete-time optimal control problem with time-varying dynamics and cost functions, where state transitions are subject to random disturbances:
\begin{align}\label{equ:dynamics-and-costs:general}
    &\text{Control dynamics: } X_{t+1} = f_t\left(X_t, U_t; W_t\right), \quad 0\leq t < T, \text{ with the initial state } X_0 = x_0;\nonumber\\
    &\text{Stage cost: } h_t(X_t, U_t), \quad 0\leq t < T, \quad \text{and terminal cost: } h_T(X_T).
\end{align}
At each time step $t$, we let $X_t$ denote the system state and $U_t$ denote the control action chosen by an agent. The function $f_t: \mathbb{R}^n \times \mathbb{R}^m \times \mathbb{R}^k \to \mathbb{R}^n$ defines how the next state $X_{t+1}$ depends on the current state $X_t$, the control action $U_t$, and the random disturbance $W_t$. The agent incurs a stage cost $h_t(X_t, U_t)$ at each time step $t < T$ and a terminal cost $h_T(X_T)$ at the final time step $T$. At each time step $t$, the controller observes the past disturbance $W_{t-1}$ and a (possibly random) prediction vector $V_t(\theta) \in \mathbb{R}^d$ before selecting a control action $U_t$, where $\theta$ is a parameter of the predictor generating the prediction. We formally define the concept of \textit{predictions} and the parameter $\theta$ in the following.

\begin{definition}[Predictions]\label{def:predictions}
At each time step $t$, the predictor with parameter $\theta \in \Theta$ provides a prediction $V_t(\theta)$, where $\Theta$ denotes the set of all possible predictor parameters. The predictions $\{V_{0:T-1}(\theta)\}_{\theta \in \Theta}$ and the disturbances $W_{0:T-1}$ live in the same probability space.
\end{definition}

Compared with previous works \citep{lin2021perturbation,li2022robustness} that assume predictions targeting specific disturbances, Definition \ref{def:predictions} focuses on the stochastic relationship between predictions and system uncertainties, yielding a unified framework for comparing different forms of prediction based on their effectiveness for control---even if their precise nature is unknown. Because predictions and disturbances share the same probability space, we can compare prediction sequences $V_{0:T-1}(\theta)$ and $V_{0:T-1}(\theta')$, generated by different predictors with parameters $\theta$ and $\theta'$.

Observe that the disturbances \(W_{0:T-1}\) and predictions in \Cref{def:predictions} do not depend on the current state or past trajectory, reflecting their exogenous nature. For example, consider the problem of quadcopter control in windy conditions \citep{o2022neural}. In this case, the wind disturbances are not influenced by the quadcopter’s state or control inputs. Under this causal relationship, we define the \emph{problem instance} as $\Xi = \bigl(W_{0:T-1}, \{V_{0:T-1}(\theta)\}_{\theta\in\Theta}\bigr)$,
and make the following assumption.

\begin{assumption}\label{assump:oblivious-disturbances-predictions}
The problem instance $\Xi$ is sampled from the distribution of problem instances before the control process starts, i.e., it will not be affected by the controller's states/actions.
\end{assumption}

Let \(\xi = \bigl(w_{0:T-1}, \{v_{0:T-1}(\theta)\}_{\theta \in \Theta}\bigr)\) denote a realization of the problem instance, including disturbances and all parameterized predictions. Under Assumption \ref{assump:oblivious-disturbances-predictions}, \(\Xi\) is viewed as realized to \(\xi\) before control begins, although the agent observes each disturbance and prediction step by step. Similar assumptions about oblivious environments or predictions appear in online optimization \citep{hazan2016introduction,rutten2023smoothed}, ensuring that future disturbances or predictions will not be affected by past states or actions. Hence, for a fixed predictor parameter \(\theta\), we define a \emph{predictive policy} as a mapping from the current state and past disturbances and predictions to a control action.

\begin{definition}[Predictive policy]\label{def:predictive-policy}
Consider a fixed predictor parameter $\theta$. For each time step $t$, let $I_t(\theta) \coloneqq (W_{0:t-1}, V_{0:t}(\theta))$ denote the history of past disturbances and predictions, and let $\mathcal{F}_t(\theta) \coloneqq \sigma(I_t(\theta))$\footnote{For any random variable $Y$, we use $\sigma(Y)$ to denote the $\sigma$-algebra it generates.}. A \emph{predictive policy} that applies to the predictor with parameter $\theta$ is a sequence of functions $\pi_{0:T-1}$, where $\pi_t$ maps a state/history pair to a control action.
\end{definition}

Given a fixed predictive policy sequence $\mathbf{\pi} = \pi_{0:T-1}$ for a predictor parameter $\theta$, we evaluate its performance via the expected total cost over $\Xi$: 
$J^\mathbf{\pi}(\theta) \coloneqq \mathbb{E}[\sum_{t=0}^{T-1} h_t(X_t, U_t) + h_T(X_T)]$, where $X_0 = x_0$, $X_{t+1} = f_t(X_t, U_t; W_t)$, $U_t = \pi_t(X_t; I_t(\theta))$, for $t=0,\ldots,T-1$. The optimal cost under $\theta$ is defined as $J^*(\theta) = \min_{\pi} J^\pi(\theta),$ where the minimum is over all predictive policies that use the predictor parameter $\theta$.

Following \citep{yu2020power}, we define \emph{prediction power} by comparing against a baseline that provides minimal information (\textit{e.g.}, no prediction). Without loss of generality, let $\mathbf{0} \in \Theta$ be the baseline predictor parameter
% \footnote{For example, to model “no prediction,” one can set $V_t(\mathbf{0}) = 0$ for all $t$.}, 
so that any $\theta \neq \mathbf{0}$ provides at least as much information as $\mathbf{0}$, \textit{i.e.}, $\mathcal{F}_t(\theta) \supseteq \mathcal{F}_t(\mathbf{0})$. Based on this baseline, we define \emph{prediction power} as the maximum possible cost improvement achieved by using predictions under $\theta$ relative to the baseline, formally stated in Definition~\ref{def:prediction-power}.

\begin{definition}[Prediction power]\label{def:prediction-power}
For a predictor with parameter $\theta$, its prediction power in the optimal control problem \eqref{equ:dynamics-and-costs:general} is $P(\theta) \coloneqq J^*(\mathbf{0}) - J^*(\theta).$
\end{definition}

Our definition of prediction power is based on the optimal control policy under a given predictor parameter and, therefore, is independent of any specific policy class.  Many previous works have considered prediction-enabled improvement within a specific policy class \citep{li2018using,li2019online,chen2015online}, where they focus on changes in $J^\pi(\theta)$ rather than $J^*(\theta)$. In other works, policies include parameters that can be tuned to perform optimally under a specific predictor; that is, $\min_{\pi \in \text{a policy class}} J^\pi(\theta)$. While these approaches are useful in specific application scenarios, our definition, based on the general optimal policy, is more universal because: (1) imposing policy class constraints may lead to performance loss, and (2) the extent of improvement can depend on policy design and parameterization, which shifts the focus away from valuing predictions themselves.

Throughout this paper, we use $\bar{\pi} = \bar{\pi}_{0:T-1}$ and $\pi^\theta = \pi_{0:T-1}^\theta$ to denote the optimal policy for the predictor with parameter $\mathbf{0}$ and $\theta$ respectively. In other words, $J^{\bar{\pi}}(\mathbf{0}) = J^*(\mathbf{0})$ and $J^{\pi^\theta}(\theta) = J^*(\theta)$. To compare the policies $\pi^\theta$ and $\bar{\pi}$, we introduce the \emph{instance-dependent Q function}, inspired by the Q function in the study of Markov decision processes (MDPs). For a given state-action pair $(x,u)$ and problem instance $\xi$, the instance-dependent Q function for a policy $\pi$ evaluates the remaining cost incurred by taking action $u$ from state $x$ and then following policy $\pi$ for all future time steps. Using $\iota_{\tau}(\theta)$ to denote the realization of $I_{\tau}(\theta)$, the instance-dependent Q function is defined as
\begin{align}\label{equ:cost-to-go-function}\textstyle
    Q_t^{\pi^\theta}(x, u; \xi) = \sum_{\tau=t}^{T-1} h_\tau(x_\tau, u_\tau) + h_T(x_T), \ \text{ where } x_t = x, \ u_t = u,
\end{align}
subject to the constraints that $x_{\tau+1} = f_\tau(x_\tau, u_\tau; w_\tau)$ for $t\leq \tau < T$ and $u_\tau = \pi_\tau^\theta\left(x_\tau; \iota_\tau(\theta)\right)$ for $t < \tau < T$. The disturbance $w_\tau$ and the history $\iota_\tau(\theta)$ in \eqref{equ:cost-to-go-function} are decided by the problem instance $\xi$, which is an input to $Q_t^{\pi^\theta}$. Similarly, we can define $Q_t^{\bar{\pi}}(x, u; \xi)$ by replacing $\theta$ with $\mathbf{0}$ and $\pi^\theta$ with $\bar{\pi}$ in \eqref{equ:cost-to-go-function}. Importantly, our instance-dependent Q function is different from the classical definition of the Q function for MDPs or reinforcement learning (RL), where it is the \textit{expectation} of the cost to go. The instance-dependent Q function denotes the actual remaining cost, which is a $\sigma(\Xi)$-measurable \textit{random variable}. The classic definition of the Q function can be recovered by taking the conditional expectation, \textit{i.e.}, $\mathbb{E}\left[Q_t^{\pi^\theta}(x, u; \Xi)\mid I_t(\theta) = \iota_t(\theta)\right]$. It is worth noting that our instance-dependent Q function is about the \textit{cost} instead of the \textit{reward}, so lower values are better.

With this definition of the instance-dependent Q function, the policies $\bar{\pi}$ and $\pi^\theta$ can be expressed as recursively minimizing the corresponding expected Q functions conditioned on the available history. Starting with $C_T^{\pi^\theta}(x; \xi) = h_T(x)$, for time step $t = T-1, \ldots, 0$, we have
\begin{align}\label{equ:Q-function-optimal-policy-recursive}
    &Q_t^{\pi^\theta}(x, u; \xi) \coloneqq h_t(x, u) + C_{t+1}^{\pi^\theta}(f_t(x, u; w_t); \xi),\text{ for } x \in \mathbb{R}^n,\ u \in \mathbb{R}^m, \text{ and problem instance }\xi;\nonumber\\
    &\pi_t^\theta(x; \iota_t(\theta)) \coloneqq \argmin_{u \in \mathbb{R}^m} \mathbb{E}\left[Q_t^{\pi^\theta}(x, u; \Xi)\mid I_t(\theta) = \iota_t(\theta)\right], \text{ for } x \in \mathbb{R}^n \text{ and history }\iota_t(\theta);\nonumber\\
    &C_t^{\pi^\theta}(x; \xi) \coloneqq Q_t^{\pi^\theta}(x, \pi_t^\theta(x; \iota_t(\theta)); \xi),  \text{ for } x \in \mathbb{R}^n \text{ and problem instance }\xi.
\end{align}
Similar recursive relationships also defines the optimal policy $\bar{\pi}$ for the baseline predictions, and we only need to replace $\theta$ with $\mathbf{0}$ and $\pi^\theta$ with $\bar{\pi}$ in the above equations. The recursive equations in (\ref{equ:Q-function-optimal-policy-recursive}) can be viewed as a generalization of the classical Bellman optimality equation for general MDPs.

\section{LTV Dynamics with Quadratic Costs}\label{sec:LQR}
We first characterize the prediction power (\Cref{def:prediction-power}) in a linear time-varying (LTV) dynamical system with quadratic costs, where the dynamics and costs are given by:
\begin{align}\label{equ:TV-LQR}
    &\text{Control dynamics: } X_{t+1} = A_t X_t + B_t U_t + W_t, \text{ for } 0\leq t < T; \nonumber\\
    &\text{stage cost: } X_t^\top Q_t X_t + U_t^\top R_t U_t, \text{ for } 0\leq t < T; \text{and terminal cost: } X_T^\top P_T X_T,
\end{align}
where $Q_{0:T-1}, R_{0:T-1},$ and $P_T$ are symmetric positive definite. The classic linear quadratic regulator (LQR) problem, along with its time-varying variant that we consider, has been used widely as a benchmark setting in the learning-for-control literature. It also serves as a good approximation of nonlinear systems near equilibrium points, making it amenable to standard analytical tools. We begin by defining key quantities that will be useful for stating the main results in this section.
For $t = T-1, \ldots, 0$, we define the matrices $H_t$, $P_t$, and $K_t$ recursively according to
\begin{align}\label{def:LTV-PKH:e1}
    H_t ={}& B_t (R_t + B_t^\top P_{t+1} B_t)^{-1} B_t^\top,\quad P_t = Q_t + A_t^\top P_{t+1} A_t - A_t^\top P_{t+1} H_t P_{t+1} A_t,\text{ and }\nonumber\\
    K_t ={}& (R_t + B_t^\top P_{t+1} B_t)^{-1} (B_t^\top P_{t+1} A_t).
\end{align}
Moreover, we define the transition matrix $\Phi_{t_2, t_1}$ as $\Phi_{t_2, t_1} = I$ if $t_2\leq t_1$ and
\begin{align}\label{def:LTV-PKH:e1-1}
    \Phi_{t_2, t_1} =
        (A_{t_2-1} - B_{t_2-1} K_{t_2-1}) (A_{t_2-2} - B_{t_2-2} K_{t_2-2}) \cdots (A_{t_1} - B_{t_1} K_{t_1}),\ \text{if } t_2 > t_1.
\end{align}
The matrix $K_t$ is the feedback gain matrix in the optimal policy, and $P_t$ is the matrix that defines the quadratic term in the optimal cost-to-go function. To simplify notation, we define the shorthands $W_{\tau\mid t}^\theta \coloneqq \mathbb{E}\left[W_\tau\mid I_t(\theta)\right]$ and $w_{\tau\mid t}^\theta \coloneqq \mathbb{E}\left[W_\tau\mid I_t(\theta) = \iota_t(\theta)\right]$.

\begin{proposition}\label{thm:LQR-closed-form}
In the case of LTV dynamics with quadratic costs, the conditional expectation of the optimal Q function $\mathbb{E}\left[Q_t^{\pi^\theta}(x, u; \Xi)\mid I_t(\theta) = \iota_t(\theta)\right]$ can be expressed as
\begin{align*}\textstyle
    \left(u + K_t x - \bar{u}_t^\theta(\iota_t(\theta))\right)^\top (R_t + B_t^\top P_{t+1} B_t) \left(u + K_t x - \bar{u}_t^\theta(\iota_t(\theta))\right) + \psi_t^{\pi^\theta}(x; \iota_t(\theta)),
\end{align*}
where $\psi_t^{\pi^\theta}(x; \iota_t(\theta))$ is a function of the state $x$ and the history $\iota_t(\theta)$ that does not depend on the control action $u$. Here, $\bar{u}_t^\theta(\iota_t(\theta)) \coloneqq - (R_t + B_t^\top P_{t+1} B_t)^{-1} B_t^\top \sum_{\tau = t}^{T-1} \Phi_{\tau+1, t+1}^\top P_{\tau + 1} w_{\tau\mid t}^\theta.$
And the optimal policy can be expressed as $\pi_t^\theta(x; \iota_t(\theta)) = - K_t x + \bar{u}_t^\theta(\iota_t(\theta)).$
\end{proposition}

We derive the closed-form expressions in \Cref{thm:LQR-closed-form} by induction following the backward recursive equations in \eqref{equ:Q-function-optimal-policy-recursive}; the full proof is deferred to \Cref{appendix:thm:LQR-closed-form}. With these expressions, we obtain a closed-form expression of the prediction power. We defer its proof to \Cref{appendix:thm:pred-power-LTV}.

\begin{theorem}\label{thm:pred-power-LTV}
In the case of LTV dynamics with quadratic costs, the prediction power of the predictor with parameter $\theta$ is $P(\theta) = \sum_{t=0}^{T-1} \Tr{(R_t + B_t^\top P_{t+1} B_t) \mathbb{E}\left[\Cov{\bar{u}_t^\theta(I_t(\theta))\mid \mathcal{F}_t(\zero)}\right]}$.
\end{theorem}

While the optimal policy in \Cref{thm:LQR-closed-form} is restricted to the LQR case, we can interpret the optimal policy as planning according the conditional expectation following the idea of model predictive control (MPC) \citep{yu2020power}, which is easier to generalize. The agent needs to solve an optimization problem and re-plan at every time step. At time step $t$, the agent solves
\begin{equation}\label{equ:MPC-in-expectation}
% \textstyle
    \argmin_{u_{t:T-1}}\quad
    \mathbb{E}\left[\sum_{\tau = t}^{T-1} h_\tau(X_\tau, u_\tau) + h_T(X_T) \given I_t(\theta) = \iota_t(\theta)\right]
\end{equation}
subject to the constraints that $X_{\tau+1} = f_\tau(X_\tau, u_\tau; W_\tau)$ for $\tau \geq t$ and $X_t = x$. Then, the agent commits to the first entry $u_{t\mid t}$ of the optimal solution as $\pi_t^\theta(x; \iota_t(\theta))$. In the LQR setting, we can further simplify it to be \textit{planning according to $w_{\tau\mid t}^\theta$} (see \Cref{appendix:MPC-form}).

The MPC forms of the optimal policy in \eqref{equ:MPC-in-expectation} extends the result in \cite{yu2020power}, which shows that MPC is the optimal predictive policy under the accurate prediction model in time-variant LQR. When the predictions are inaccurate, and the system is time-varying, MPC is still optimal if we solve the predictive optimal control problem in expectation \eqref{equ:MPC-in-expectation}.

\textbf{Evaluation.} One can follow the expressions in \Cref{thm:pred-power-LTV} to evaluate the prediction power, but it requires taking the conditional covariance on the top of conditional expectations ($\bar{u}_t^\theta(\iota_t(\theta))$ in \Cref{thm:LQR-closed-form}). To avoid this recursive structure, an alternative way is to first construct the \textit{surrogate optimal action}, which is defined as
\begin{align}\label{equ:surrogate-optimal-action}\textstyle
    \bar{u}_t^*(\Xi) \coloneqq - (R_t + B_t^\top P_{t+1} B_t)^{-1} B_t^\top \sum_{\tau = t}^{T-1} \Phi_{\tau+1, t+1}^\top P_{\tau + 1} W_\tau.
\end{align}
We call $\bar{u}_t^*(\Xi)$ the surrogate-optimal action, because it is the optimal action that an agent should take with the oracle knowledge of all future disturbances at time $t$. The prediction power in \Cref{thm:pred-power-LTV} can be expressed as $\mathbb{E}\left[\Cov{\bar{u}_t^{*}(\Xi)\mid I_t(\mathbf{0})}\right] - \mathbb{E}\left[\Cov{\bar{u}_t^{*}(\Xi)\mid I_t(\theta)}\right]$. Following this decomposition, we propose an evaluation approach that constructs $\bar{u}_t^*(\Xi)$ before estimating its conditional covariance with respect to $I_t(\theta)$ and $I_t(\mathbf{0})$ separately. We defer the details to \Cref{appendix:alg:pred-power-eva}.

\subsection{Prediction Power \texorpdfstring{$\neq$}{≠} Accuracy}\label{sec:main:power-vs-accuracy}
As \Cref{thm:LQR-closed-form} suggests, one way to implement the optimal policy is to predict each of the future disturbances $W_{t:T-1}$ and generate the estimations $w_{(t:T-1)\mid t}^\theta$ in deciding the action at time step $t$. However, two controllers with the same estimation error (as measured by mean squared error (MSE)) can have very different control costs. Because of this reason, the control cost bounds depend on the estimation errors in previous works \citep{zhang2021regret,yu2022competitive,lin2022bounded} must be loose, so one cannot rely on them to infer or compare the values of different predictors.

To illustrate this point, we provide an example where the prediction power can change significantly when the prediction accuracy does not change.

\begin{example}\label{example:one-step-dim-rotation}
Consider the time-invariant LQR setting, \emph{i.e.}, assume $A_t = A, B_t = B, Q_t = Q, R_t = R$ for all $t$ and $P_T = P$ is the solution to the Discrete-time Riccati Equation (DARE) in \eqref{equ:TV-LQR}. Suppose the disturbance is sampled $W_t \overset{\text{i.i.d.}}{\sim} N(0, I)$ at every time step $t$. Let $\rho \in [0, \frac{\sqrt{2}}{2}]$ be a fixed coefficient. We construct a class of predictors from the disturbances $\{W_t\}$ by applying the affine transformation $V_t(\theta) \coloneq \rho \theta W_t + \epsilon_t(\rho, \theta)$ for $\theta \in \mathbb{R}^{2\times 2}$ that satisfies $\theta \theta^\top \preceq \frac{1}{2}I$, where the random noise $\epsilon_t(\rho, \theta)$ is independently sampled from a Gaussian distribution $N(0, I - \rho^2 \theta \theta^\top)$.

We can construct $\theta$ such that $V_t(\theta)$ and $V_t(I)$ achieve the same mean-square error (MSE) when predicting each individual entry of $W_t$, yet $P(I) > P(\theta)$. To construct $\theta$, note that $(W_t, V_t(\theta))$ satisfies $\mathbb{E}\left[W_t\mid V_t(\theta)\right] = \rho \theta^\top V_t$ and $\Cov{W_t\mid V_t(\theta)} = I - \rho^2 \theta^\top \theta$. Thus, we can change $\theta$ without affecting the MSE of predicting each individual entry as long as the diagonal entries of $\theta^\top \theta$ remain the same. However, by \Cref{thm:pred-power-LTV}, we know the prediction power is equal to $\rho^2 T \cdot \Tr{\theta^\top \theta P H P}$, where $H = B (R + B^\top P B)^{-1} B^\top$. Thus, the off-diagonal entries of $\theta^\top \theta$ can also affect the value of $\Tr{\theta^\top \theta P H P}$. We instantiate this example with a 2-D double-integrator dynamical system in Appendix \ref{appendix:example-details:double-integrator}: the predictors with parameters $I$ and $\theta$ shares the same MSE but their prediction powers are significantly different.\footnote{The simulation code for all examples (Examples \ref{example:one-step-dim-rotation}, \ref{example:pred-power-and-policy-opt}, and \ref{example:multi-step-1d}) can be found at \url{https://github.com/yihenglin97/Prediction-Power}.}
\end{example}

Example \ref{example:one-step-dim-rotation} shows how prediction power can vary even when the accuracy of predicting each entry of the disturbance $W_t$ remains the same, where the construction leverages the covariance between the predictions for different entries of $W_t$. While the construction in Example \ref{example:one-step-dim-rotation} requires $n \geq 2$, we also provide an example with $n = 1$ and multiple steps of predictions in Appendix \ref{appendix:example-details}. From these examples, it is clear that one should not use the MSEs of predicting future disturbances to infer the prediction power. The intuition behind this mismatch is that MSE does not depend on matrices $(A, B, Q, R)$, but the prediction power does. The mismatch also relates to the findings in the decision-focused learning literature discussed in the related work section.

\subsection{Prediction Power and Online Policy Optimization}\label{sec:main:power-and-optimization}
The closed-form expression of the prediction power, presented in \Cref{thm:pred-power-LTV}, characterizes the maximum potential of using a given prediction sequence $V_{0:T-1}(\theta)$. Here, we draw a connection between prediction power and online policy optimization \cite{lin2023online,agarwal2019online}, which aims to learn and adapt the optimal control policy within a certain policy class over time: the prediction power serves as an improvement upper bound of applying online policy optimization to predictive policies, although it is generally unattainable. In the following example, we demonstrate this bound using M-GAPS \cite{lin2024online}, a state-of-the-art online policy optimization algorithm.

\begin{figure*}
\centering
\begin{minipage}{.48\textwidth}
  \centering
  \includegraphics[width=.9\linewidth]{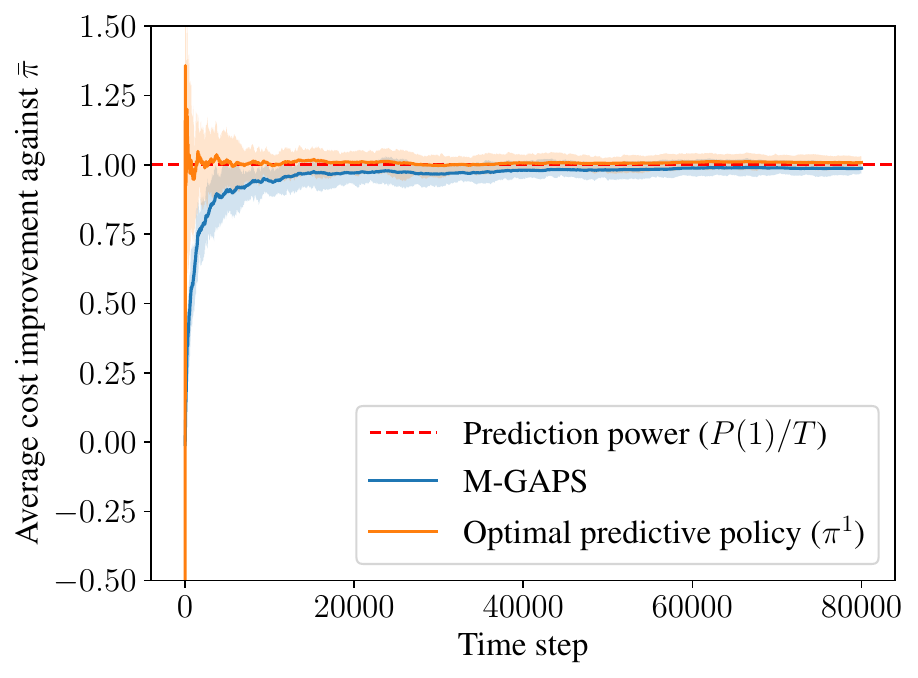}
  \caption{\Cref{example:pred-power-and-policy-opt}: Prediction $V_t(1)$ is available. Candidate policy: $u_t = -Kx_t + \Upsilon_t v_t(1)$.}
  \label{fig:Pred-power-vs-avg-cost:Current}
\end{minipage}%
\hfill
\begin{minipage}{.48\textwidth}
  \centering
  \includegraphics[width=.9\linewidth]{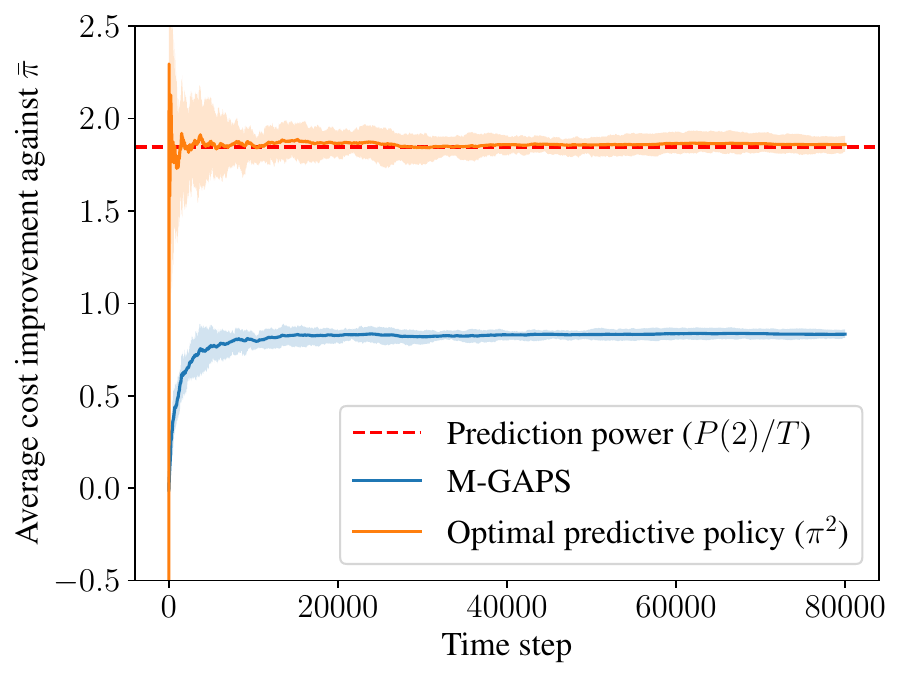}
  \caption{\Cref{example:pred-power-and-policy-opt}: Prediction $V_t(2)$ is available. Candidate policy: $u_t = -Kx_t + \Upsilon_t v_t(2)$.}
  \label{fig:Pred-power-vs-avg-cost:Next}
\end{minipage}
\end{figure*}

\begin{example}\label{example:pred-power-and-policy-opt}
We construct two scenarios under the same setting as \Cref{example:one-step-dim-rotation}. First, when the prediction is $V_t(1)\coloneqq \rho W_t + \epsilon_t(\rho, I)$, we let M-GAPS adapt within the candidate policy class $u_t = - K x_t + \Upsilon_t v_t(1)$, where $\Upsilon_t \in \mathbb{R}^{1\times 2}$ is the policy parameter. Here, the optimal predictive policy $\pi^1$ is contained in the candidate policy class. We plot the average cost improvement of M-GAPS and $\pi^1$ compared against the optimal no-prediction policy $\bar{\pi}$ in \Cref{fig:Pred-power-vs-avg-cost:Current}. From the initialization $\Upsilon_0 = \mathbf{0}$, M-GAPS tunes $\Upsilon_t$ to improve the average cost over time, and the average cost improvement against $\bar{\pi}$ converges towards the averaged prediction power $P(1)/T$.

In the second scenario, we change the prediction to apply M-GAPS to $V_t(2)\coloneqq V_{t+1}(1)$ (i.e., the same prediction as before is made available 1-step ahead). We let M-GAPS adapts within the same candidate policy class $u_t = - K x_t + \Upsilon_t v_t(2)$, where the policy parameter is still $\Upsilon_t \in \mathbb{R}^{1\times 2}$. Unlike the first scenario, the optimal predictive policy $\pi^2$ is \emph{not} contained in the candidate policy class, because $\pi^2$ uses both $v_t(2)$ and $v_{t-1}(2)$ to decide the action. As a result, M-GAPS cannot achieve an improvement that is close to the averaged prediction power $P(2)/T$, which is achievable by $\pi^2$ (see \Cref{fig:Pred-power-vs-avg-cost:Next}).
\end{example}

The details of \Cref{example:pred-power-and-policy-opt} are provided in \Cref{appendix:example-pred-power-and-policy-opt}. It demonstrates how prediction power serves as an upper bound for the cost improvement achieved by online policy optimization. Conversely, online policy optimization offers practical tools to achieve (part of) the potential benefit of using predictions without requiring explicit knowledge or estimation of the joint distribution between predictions and true disturbances.

\section{Characterizing Prediction Power: A General-Purpose Theorem}\label{sec:main_results}
In this section, we provide a theorem to characterize the prediction power $P(\theta)$ within the general problem setting introduced in \Cref{sec:setting}. Our result relies on two conditions about a growth property of the expected Q function under $\pi^\theta$ and the covariance of the optimal policy's action when conditioned on the $\sigma$-algebra $\mathcal{F}_t(\mathbf{0})$ of the baseline. We state these conditions and provide intuitive explanations.

\begin{condition}\label{cond:Q-function-quadratic-growth}
For a sequence of positive semi-definite matrices $M_{0:T-1}$, the following inequality holds for all time steps $0 \leq t < T$: For any $x \in \mathbb{R}^n$, $u \in \mathbb{R}^m$, and history $\iota_t(\theta)$,
\begin{align}\label{thm:Bellman-improvement-not-worse:assump-1}
    \mathbb{E}\left[Q_t^{\pi^\theta}(x, u; \Xi) - C_t^{\pi^\theta}(x; \Xi)\mid I_t(\theta) = \iota_t(\theta)\right] \geq (u - \pi_t^\theta(x; \iota_t(\theta)))^\top M_t (u - \pi_t^\theta(x; \iota_t(\theta))).
\end{align}
\end{condition}
The LQR setting (\Cref{sec:LQR}) satisfies \Cref{cond:Q-function-quadratic-growth} with $M_t = R_t + B_t^\top P_{t+1} B_t$.

Condition \ref{cond:Q-function-quadratic-growth} states that conditioned on any history $\iota_t(\theta)$, the expected Q function of policy $\pi^\theta$ grows at least quadratically as the action $u$ deviates from the optimal policy's action. Note that one can always pick $M_t$ to be the all-zeros matrix to make Condition \ref{cond:Q-function-quadratic-growth} hold, but the choice of $M_t$ will affect the prediction power bound in Theorem \ref{thm:Bellman-improvement-not-worse}. When $M_t \succ 0$, deviating from the action of policy $\pi^\theta$ causes a non-negligible loss. The loss is characterized by the difference between the resulting Q function value and the cost-to-go function value. When this condition does not hold with any non-zero matrix $M_t$, one can construct an extreme case when $Q_t^{\pi^\theta}$ is a constant by letting all cost functions $h_{0:T}$ be constants; in this case, the prediction power must be zero because every policy achieves the same total cost no matter what predictions they use.

\begin{condition}\label{cond:optimal-policy-psd-covariance}
One of the following holds for the optimal policy $\pi^\theta$:
\begin{enumerate}[nolistsep, leftmargin=*, label=(\alph*)]
\item For positive semi-definite matrices $\Sigma_{0:T-1}$, the following holds for all time steps $0 \leq t < T$: \label{cond:optimal-policy-psd-covariance:a}
\begin{align}\label{thm:Bellman-improvement-not-worse:assump-2}
    &\mathbb{E}\left[\Cov{\pi_t^\theta(X; I_t(\theta))\mid I_t(\mathbf{0})}\right] \succeq \Sigma_t, \text{ for any }\mathcal{F}_t(\mathbf{0})\text{-measurable } X.
\end{align}
\item For nonnegative scalars $\sigma_{0:T-1}$, the following holds for all time steps $0 \leq t < T$: \label{cond:optimal-policy-psd-covariance:b}
\begin{align}\label{thm:Bellman-improvement-not-worse:assump-2-weaker}
    \mathbb{E}\left[\Tr{\Cov{\pi_t^\theta(X; I_t(\theta))\mid I_t(\mathbf{0})}}\right] \geq \sigma_t, \text{ for any }\mathcal{F}_t(\mathbf{0})\text{-measurable } X.
\end{align}
\end{enumerate}
\end{condition}

Before discussing the details, we note that by setting \(\sigma_t = \text{Tr}(\Sigma_t)\), Condition \ref{cond:optimal-policy-psd-covariance}\,\ref{cond:optimal-policy-psd-covariance:a} implies (and is therefore stronger than) Condition \ref{cond:optimal-policy-psd-covariance}\,\ref{cond:optimal-policy-psd-covariance:b}. Similar to \Cref{cond:Q-function-quadratic-growth}, one can always pick $\Sigma_t$ to be all-zeros matrix to satisfy \Cref{cond:optimal-policy-psd-covariance}\,\ref{cond:optimal-policy-psd-covariance:a}, but it will affect the prediction power bound. The LQR setting (\Cref{sec:LQR}) satisfies \Cref{cond:optimal-policy-psd-covariance}\,\ref{cond:optimal-policy-psd-covariance:a} with $\Sigma_t = \mathbb{E}\left[\Cov{\bar{u}_t^\theta(I_t(\theta))\mid \mathcal{F}_t(\zero)}\right].$

Condition \ref{cond:optimal-policy-psd-covariance}\,\ref{cond:optimal-policy-psd-covariance:a} states that conditioned on the history $I_t(\mathbf{0})$ from the baseline, the covariance matrix of policy $\pi^\theta$'s action from any $\mathcal{F}_t(\mathbf{0})$-measurable state is positive semi-definite in expectation. Recall that $\mathcal{F}_t(\mathbf{0}) = \sigma(I_t(\mathbf{0}))$. To understand this, suppose that the agent only has access to the baseline information. Then, the agent cannot predict the action that policy $\pi^\theta$ would take. This should usually hold because the action $\pi_t^\theta(X; I_t(\theta))$ is not $\mathcal{F}_t(\mathbf{0})$-measurable, and the lower bound in \eqref{thm:Bellman-improvement-not-worse:assump-2} implies the mean-square prediction error cannot improve below a certain threshold. When this condition does not hold with non-zero matrix $\Sigma_t$ (or scalar $\sigma_t$), one can design a policy $\bar{\pi}'$ that always picks the same action as $\pi^\theta$ but only requires access to the baseline information $I_t(\mathbf{0})$, which implies $P(\theta) = 0$ because $J^*(\mathbf{0}) \leq J^{\bar{\pi}'}(\mathbf{0}) = J^*(\theta)$. This can happen, for example, when $W_{0:T-1}$ are deterministic. 

Note it is possible that the optimal action at different states has a positive variance in different directions, but there is no non-trivial lower bound on the covariance matrix as required by Condition \ref{cond:optimal-policy-psd-covariance}\,\ref{cond:optimal-policy-psd-covariance:a}. In this case, Condition \ref{cond:optimal-policy-psd-covariance}\,\ref{cond:optimal-policy-psd-covariance:b} provides a weaker alternative and would be useful when we can only establish a lower bound on the trace of the optimal action's covariance matrix (\textit{e.g.}, Section \ref{sec:LTV-well-conditioned}). 

\begin{theorem}\label{thm:Bellman-improvement-not-worse}
If Conditions \ref{cond:Q-function-quadratic-growth} and \ref{cond:optimal-policy-psd-covariance}\,\ref{cond:optimal-policy-psd-covariance:a} hold with matrices $M_{0:T-1}$ and $\Sigma_{0:T-1}$, then $P(\theta) \geq \sum_{t=0}^{T-1} \Tr{M_t \Sigma_t}$. Alternatively, if Conditions \ref{cond:Q-function-quadratic-growth} and \ref{cond:optimal-policy-psd-covariance}\,\ref{cond:optimal-policy-psd-covariance:b} hold with matrices $M_{0:T-1}$ and scalars $\sigma_{0:T-1}$, then $P(\theta) \geq \sum_{t=0}^{T-1} \mu_\text{min}(M_t) \cdot \sigma_t$, where $\mu_\text{min}(\cdot)$ returns the smallest eigenvalue.
\end{theorem}

\begin{figure}
\centering
\includegraphics[width=0.7\textwidth]{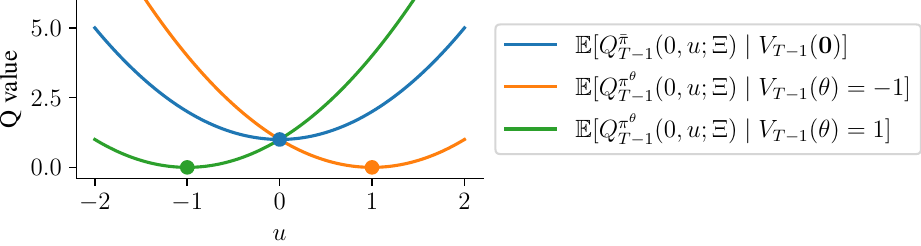}%
\vspace{-0.75em}%
\caption{An illustration of why predictions are helpful, corresponding to \Cref{ex:prediction-improvement-intuition}. The expected Q functions with perfect predictions (green and orange lines) have lower minima than the expected Q function with uninformative predictions (blue line).}
\label{fig:pred-improve-intuition}
\end{figure}

We defer the proof of Theorem \ref{thm:Bellman-improvement-not-worse} to Appendix \ref{appendix:thm:Bellman-improvement-not-worse}. As a remark, in the LQR setting, the first inequality in \Cref{thm:Bellman-improvement-not-worse} holds with equality, and it recovers the same expression as \Cref{thm:pred-power-LTV} in \Cref{sec:LQR}. There are two main takeaways of Theorem \ref{thm:Bellman-improvement-not-worse}. First, recall that one can always pick $M_t$ and $\Sigma_t$ to be the all-zeros matrices to satisfy Conditions \ref{cond:Q-function-quadratic-growth} and \ref{cond:optimal-policy-psd-covariance}. In this case, Theorem \ref{thm:Bellman-improvement-not-worse} states that $P(\theta)\geq 0$, which means that having predictions, no matter how weak they are, does not hurt. Second, to characterize the improvement in having predictions, Conditions \ref{cond:Q-function-quadratic-growth} and \ref{cond:optimal-policy-psd-covariance} can establish a lower bound for the prediction power that is strictly positive if $\Tr{M_t \Sigma_t} > 0$ or $\mu_\text{min}(M_t) \sigma_t > 0$. We provide an example to help illustrate how Conditions \ref{cond:Q-function-quadratic-growth} and \ref{cond:optimal-policy-psd-covariance}\,\ref{cond:optimal-policy-psd-covariance:a} can work together to ensure that the predictions can lead to a strict improvement on the control cost (see Figure \ref{fig:pred-improve-intuition} for an illustration). 

\begin{example}
\label{ex:prediction-improvement-intuition}
Consider the following optimal control problem
\begin{align*}
    \text{Dynamics: } X_{t+1} = U_t + W_t,\
    \text{Stage cost: } h_t(x, u) = x^2,\
    \text{Terminal cost: } h_T(x) = x^2,
\end{align*}
where each disturbance $W_t$ is sampled independently according to $\mathbb{P}(W_t = -1) = \mathbb{P}(W_t = 1) = \frac{1}{2}$. Suppose that the predictor with parameter $\theta$ can predict $W_t$ exactly (\emph{i.e.}, $V_t(\theta) = W_t$), while the baseline predictor is uninformative (\emph{e.g.}, $V_t(\zero) = 0$). The Q functions, cumulative cost, and optimal actions under each predictor are
\begin{align*}
    Q_t^{\pi^\theta}(x, u; \Xi) &= x^2 + (u + V_t(\theta))^2, &
    Q_t^{\bar\pi}(x, u; \Xi) &= x^2 + (u + W_t)^2 + (T-t-1), \\
    C_t^{\pi^\theta}(x; \Xi) &= x^2, &
    C_t^{\bar\pi}(x; \Xi) &= x^2 + (T-t), \\
    \pi_t^\theta(x; I_t(\theta)) &= -V_t(\theta) = -W_t, &
    \bar\pi_t(x; I_t(\zero)) &= 0.
\end{align*}
The Q function $Q_t^{\pi^\theta}$ is strongly convex in $u$, with \Cref{cond:Q-function-quadratic-growth} holding for any $M_t \in [0,1]$. Furthermore, the optimal action has positive variance, with \Cref{cond:optimal-policy-psd-covariance}\,\ref{cond:optimal-policy-psd-covariance:a} holding for any $\Sigma_t \in [0,1]$. Thus, by \Cref{thm:Bellman-improvement-not-worse}, the prediction power satisfies $P(\theta) \geq T$. Indeed, by comparing the cumulative cost functions, we see that the predictor with parameter $\theta$ incurs a lower cumulative cost by exactly $T$ (as expected by \Cref{thm:pred-power-LTV}).

\Cref{fig:pred-improve-intuition} illustrates the expected Q functions at time $t=T-1$ and $x=0$, which the policies $\pi_t^\theta(x;I_t(\theta))$ and $\bar\pi_t(x;I_t(\zero))$ seek to minimize. The expected Q functions with perfect predictions have lower minima than the expected Q function with uninformative predictions.
\end{example}

\Cref{thm:Bellman-improvement-not-worse} provides a useful tool to characterize the prediction power by reducing the problem of comparing two policies $\pi^\theta$ and $\bar{\pi}$ over the whole horizon to studying the properties of one policy $\pi^\theta$ at each time step. Our proof of \Cref{thm:Bellman-improvement-not-worse} follows the same intuition as the widely-used performance difference lemma in RL (see Lemma 6.1 in \cite{kakade2002approximately}), comparing the per-step ``advantage'' of $\pi^\theta$ along the trajectory of $\bar{\pi}$. When only the baseline information is available, the agent must pick a suboptimal action \eqref{thm:Bellman-improvement-not-worse:assump-2} and incur a loss \eqref{thm:Bellman-improvement-not-worse:assump-1} at each step, which accumulates to the total cost difference.

While Theorem \ref{thm:Bellman-improvement-not-worse} applies to the general dynamical system and cost functions in \eqref{equ:dynamics-and-costs:general}, the two conditions with their key coefficients $M_t$ and $\Sigma_t$ (or $\sigma_t$) still depend on the optimal Q function and the optimal policy that are implicitly defined through the recursive equations \eqref{equ:Q-function-optimal-policy-recursive}. To instantiate Theorem \ref{thm:Bellman-improvement-not-worse}, we need to derive explicit expressions of $M_t$ and $\Sigma_t$ under more specific dynamics/costs.

\subsection{LTV Dynamics with General Costs}\label{sec:LTV-well-conditioned}
In this section, we consider an online optimal control problem with linear time-varying dynamics and more general cost functions compared with the LQR setting in \Cref{sec:LQR}.
\begin{align}\label{equ:setting-LTV-well-conditioned}
    &\text{Control dynamics: } X_{t+1} = A_t X_t + B_t U_t + W_t, \text{ for } 0\leq t < T;\nonumber\\
    &\text{stage cost: } h_t^x(X_t) + h_t^u(U_t), \text{ for } 0\leq t < T;
    \text{and terminal cost: } h_T^x(X_T).
\end{align}
The LTV system with quadratic cost functions studied in Section \ref{sec:LQR} is a special case of \eqref{equ:setting-LTV-well-conditioned}. The setting is challenging because the optimal Q function/policy $\pi^\theta$ do not have closed-form expressions like Proposition \ref{thm:LQR-closed-form}. To tackle it, we follow the recursive equations \eqref{equ:Q-function-optimal-policy-recursive} to establish Conditions \ref{cond:Q-function-quadratic-growth} and \ref{cond:optimal-policy-psd-covariance}\,\ref{cond:optimal-policy-psd-covariance:b}. We make the following assumptions about the cost functions and dynamical matrices:
\begin{assumption}\label{assump:well-conditioned-cost}
For every time step $t$, $h_t^x$ is $\mu_x$-strongly convex and $\ell_x$-smooth; $h_t^u$ is $\mu_u$-strongly convex and $\ell_u$-smooth; The dynamical matrices satisfy that $\mu_A I \preceq A_t^\top A_t \preceq \ell_A I$ and $\mu_B I \preceq B_t^\top B_t \preceq \ell_B I$. Further, we assume $\ell_A < 1$.
\end{assumption}

Under \Cref{assump:well-conditioned-cost}, we can verify \Cref{cond:Q-function-quadratic-growth} and show that the expected cost-to-go functions are well-conditioned. We state this result in \Cref{lemma:well-condition-assump-verification-strongly-convex} and defer its proof to \Cref{appendix:lemma:well-condition-assump-verification-strongly-convex}.

\begin{lemma}\label{lemma:well-condition-assump-verification-strongly-convex}
Under Assumption \ref{assump:well-conditioned-cost}, Condition \ref{cond:Q-function-quadratic-growth} holds with $M_t = \mu_u I$. Further,  conditional expectation $\mathbb{E}[C_t^{\pi^\theta}(x; \Xi)\mid I_t(\theta) = \iota_t(\theta)]$ as a function of $x$ is $\mu_t$-strongly convex and $\ell_t$-smooth for any history $\iota_t(\theta)$, where $\mu_t$ and $\ell_t$ are defined as following: Let $\mu_T = \mu_x$ and $\ell_T = \ell_x$,
\begin{align}\label{lemma:well-condition-assump-verification-strongly-convex:recursive-def}
    \mu_t = \mu_x + \mu_A \cdot \frac{\mu_u \mu_{t+1}}{\mu_u + b^2 \mu_{t+1}}, \text{ and }
    \ell_t = \ell_x + \ell_A \cdot \ell_{t+1}, \text{ for time } t = T-1, \ldots, 0.
\end{align}
\end{lemma}

To establish the second condition about the covariance of the optimal policy's action, we make the following assumption about the joint distribution of the disturbances and the predictions:
\begin{assumption}\label{assump:one-step-prediction-indep-pairs}
The disturbances and predictions can be grouped as pairs $\{(W_t, V_t(\theta))\}_{t=0}^{T-1}$, where $(W_t, V_t(\theta))$ is joint Gaussian and independent with $(W_{t'}, V_{t'}(\theta))$ when $t \not = t'$. Further, assume that the baseline is no prediction, i.e., $V_t(\mathbf{0}) = 0$. And for $\theta \in \Theta$, there exists $\lambda_t(\theta) \in \mathbb{R}_{\geq 0}$ such that $\Cov{W_t} - \Cov{W_t\mid V_t(\theta)}  \succeq \lambda_t(\theta) I, $ for any $0\leq t < T.$
\end{assumption}

With \Cref{assump:one-step-prediction-indep-pairs} and \Cref{lemma:well-condition-assump-verification-strongly-convex}, we can verify that Condition \ref{cond:optimal-policy-psd-covariance}\,\ref{cond:optimal-policy-psd-covariance:b} holds with
\begin{align}\label{equ:one-step-pred-second-cond}
    \Tr{\Cov{\pi_t^\theta(x; I_t(\theta))\mid \mathcal{F}_t(0)}} \geq \sigma_t \coloneqq \frac{n \lambda_t(\theta) \mu_{t+1}^2 \cdot \mu_B}{2(\ell_u + \ell_{t+1}\sqrt{\ell_B})^2}.
\end{align}
Since Conditions \ref{cond:Q-function-quadratic-growth} and \ref{cond:optimal-policy-psd-covariance}\,\ref{cond:optimal-policy-psd-covariance:b} hold, we apply \Cref{thm:Bellman-improvement-not-worse} to obtain the prediction power bound.

\begin{theorem}\label{coro:one-step-pred-power-lower}
In the case of LTV dynamics with well-conditioned costs, suppose Assumptions \ref{assump:well-conditioned-cost} and \ref{assump:one-step-prediction-indep-pairs} hold. The prediction power of the predictor with parameter $\theta$ is lower bounded by $P(\theta) \geq \sum_{t=0}^{T-1} \mu_u \sigma_t$, where $\sigma_t$ is defined in \eqref{equ:one-step-pred-second-cond}.
\end{theorem}

We provide a more detailed proof outline and the proofs in \Cref{appendix:LTV-general-costs-outline}. As a remark, the lower bound of the prediction power in Theorem \ref{coro:one-step-pred-power-lower} shows that even weak predictions (\textit{i.e.}, small but non-zero $\lambda_t(\theta)$ in Assumption \ref{assump:one-step-prediction-indep-pairs}) can help improve the control cost compared with the no-prediction baseline. Although Assumption \ref{assump:one-step-prediction-indep-pairs} limits $V_t(\theta)$ to be only correlated with $W_t$, we provide a roadmap towards more general dependencies on all future $W_{t:T-1}$ in Appendix \ref{appendix:multi-step-prediction-well-cond-roadmap}.

\section{Concluding Remarks}\label{sec:Concluding-remarks}
In this work, we propose the metric of prediction power and characterize it in the time-varying LQR setting (Theorem \ref{thm:pred-power-LTV}). We extend our analysis to provide a lower bound for the general setting (\Cref{thm:Bellman-improvement-not-worse}), which is helpful for establishing the incremental value of (weak) predictions beyond LQR (Theorem \ref{coro:one-step-pred-power-lower}). We emphasize that our framework is very broad. For example, if we let the parameter $\theta$ represent the dataset that the predictor is trained on, then the prediction power $P(\theta)$ effectively quantifies the value of that particular dataset with respect to the optimal control problem.

We would like to highlight two directions of research inspired by our results. First, while our work establishes prediction power of a predictor with parameter $\theta$ relative to a strictly less-informative baseline, it does not immediately enable comparison between two arbitrary parameters $\theta$ and $\theta'$ when our general lower bounds in \Cref{sec:main_results} are not tight. Second, while we discuss about how to evaluate prediction power of a given parameter $\theta$, our work does not specify what the optimal $\theta$ is. The problem of learning the parameter that maximizes $P(\theta)$ may be interesting future work.

\bibliographystyle{apalike}
\bibliography{main}

\newpage

\begin{center}
    {\LARGE\bfseries Appendices}
\end{center}

\appendix

\section{Proofs and Examples for LTV Dynamics with Quadratic Costs}
\subsection{Proof of Proposition \ref{thm:LQR-closed-form}}\label{appendix:thm:LQR-closed-form}
Recall that we introduce the shorthand
\begin{align*}
    W_{\tau\mid t}^\theta = \Condexp{W_\tau}{I_t(\theta)}.
\end{align*}
We show by induction that
\begin{align*}
    &\mathbb{E}\left[Q_t^{\pi^\theta}(x, u; \Xi)\mid I_t(\theta)\right]\\
    &= \left(u + K_t x - \bar{u}_t^\theta(I_t(\theta))\right)^\top (R_t + B_t^\top P_{t+1} B_t) \left(u + K_t x - \bar{u}_t^\theta(I_t(\theta))\right) + \psi_t^{\pi^\theta}(x; I_t(\theta)),\\
    &\pi_t^\theta(x; I_t(\theta)) = - K_t x + \bar{u}_t^\theta(I_t(\theta)),
\end{align*}
together with the expression of the optimal cost-to-go function
\begin{align}\label{thm:LQR_closed_form:e1}
    \mathbb{E}\left[C_t^{\pi^\theta}(x; \Xi)\mid I_t(\theta)\right] = x^\top P_t x + 2 \left(\sum_{\tau = t}^{T-1}\Phi_{\tau+1, t}^\top P_{\tau+1} W_{\tau\mid t}^\theta\right)^\top x + \Psi_t(I_t(\theta)),
\end{align}
where recall that for $t_2 > t_1$,
\begin{align*}
    \Phi_{t_2, t_1}^\top \coloneq{}& (A_{t_1} - B_{t_1} K_{t_1})^\top \cdots (A_{t_2-1} - B_{t_2-1} K_{t_2-1})^\top \\
    ={}& (A_{t_1}^\top - A_{t_1}^\top P_{t_1 + 1} H_{t_1}) \cdots (A_{t_2 - 1}^\top - A_{t_2 - 1}^\top P_{t_2} H_{t_2 - 1}).
\end{align*}
and $\Psi_t(I_t(\theta))$ is a function of the history observations/predictions which does not depend on $x$. Note that \eqref{thm:LQR_closed_form:e1} holds when $t = T$ because $C_T^{\pi^\theta}(x; \Xi) = x^\top P_T x$.

Suppose that \eqref{thm:LQR_closed_form:e1} holds for $t+1$. Then, we have
\begin{align*}
    &\mathbb{E}\left[C_{t+1}^{\pi^\theta}(x + W_t; \Xi)\mid I_t(\theta)\right]\\
    ={}& \mathbb{E}\left[\mathbb{E}\left[C_{t+1}^{\pi^\theta}(x + W_t; \Xi)\mid I_{t+1}(\theta)\right]\mid I_t(\theta)\right]\\
    ={}& \mathbb{E}\left[\left.\left(x + W_t\right)^\top P_{t+1} \left(x + W_t\right)\right| I_t(\theta)\right] + 2\mathbb{E}\left[\left.\sum_{\tau = t+1}^{T-1} \Phi_{\tau, t+1}^\top P_{\tau+1} W_{\tau\mid t+1}^\theta\right| I_t(\theta)\right]^\top x\\
    &+ 2 \mathbb{E}\left[\left.\left(\sum_{\tau = t+1}^{T-1} \Phi_{\tau, t+1}^\top P_{\tau+1} W_{\tau\mid t+1}^\theta\right)^\top W_t\right| I_t(\theta)\right] + \mathbb{E}\left[\Psi_{t+1}(I_{t+1}(\theta))\mid I_t(\theta)\right]\\
    ={}& x^\top P_{t+1} x + 2\left(P_{t+1} W_{t\mid t}^\theta + \sum_{\tau = t+1}^{T-1} \Phi_{\tau, t+1}^\top P_{\tau+1} W_{\tau\mid t}^\theta\right)^\top x + \Tr{P_{t+1} \cdot \Cov{W_t\mid I_t(\theta)}}\\
    &+ 2 \mathbb{E}\left[\left.\left(\sum_{\tau = t+1}^{T-1} \Phi_{\tau, t+1}^\top P_{\tau + 1} W_{\tau\mid t+1}^\theta\right)^\top W_t\right| I_t(\theta)\right] + \mathbb{E}\left[\Psi_{t+1}(I_{t+1}(\theta))\mid I_t(\theta)\right].
\end{align*}
To simplify the notation, let
\begin{align*}
    \bar{\psi}_{t+1}(I_t(\theta)) \coloneqq{}& \Tr{P_{t+1} \cdot \Cov{W_t\mid I_t(\theta)}} + 2 \mathbb{E}\left[\left(\sum_{\tau = t+1}^{T-1} \Phi_{\tau, t+1}^\top P_{\tau + 1} W_{\tau\mid t+1}^\theta\right)^\top W_t \given I_t(\theta)\right]\\
    &+ \mathbb{E}\left[\Psi_{t+1}(I_{t+1}(\theta))\mid I_t(\theta)\right].
\end{align*}
We see that the expected Q function is given by
\begin{align*}
    &\mathbb{E}\left[Q_{t}^{\pi^\theta}(x, u; \Xi)\mid I_t(\theta)\right]\\
    &= x^\top Q_t x + u^\top R_t u + \mathbb{E}\left[C_{t+1}^{\pi^\theta}(A_t x + B_t u + W_t; \Xi)\mid I_t(\theta)\right]\\
    &=x^\top Q_t x + u^\top R_t u + (A_t x + B_t u)^\top P_{t+1} (A_t x + B_t u)\\
    &\quad+ 2\left(P_{t+1} W_{t\mid t}^\theta + \sum_{\tau = t+1}^{T-1} \Phi_{\tau, t+1}^\top P_{\tau+1} W_{\tau\mid t}^\theta\right)^\top (A_t x + B_t u) + \bar{\psi}_{t+1}(I_t(\theta))\\
    &= u^\top (R_t + B_t^\top P_{t+1} B_t) u + 2\left(P_{t+1} A_t x + P_{t+1} W_{t\mid t}^\theta + \sum_{\tau = t+1}^{T-1} \Phi_{\tau, t+1}^\top P_{\tau+1} W_{\tau\mid t}^\theta\right)^\top B_t u\\
    &\quad+ x^\top (Q_t + A_t^\top P_{t+1} A_t) x + 2\left(P_{t+1} W_{t\mid t}^\theta + \sum_{\tau = t+1}^{T-1} \Phi_{\tau, t+1}^\top P_{\tau+1} W_{\tau\mid t}^\theta\right)^\top A_t x + \bar{\psi}_{t+1}(I_t(\theta))\\
    &= \left(u + K_t x - \bar{u}_t^\theta(I_t(\theta))\right)^\top (R_t + B_t^\top P_{t+1} B_t) \left(u + K_t x - \bar{u}_t^\theta(I_t(\theta))\right) + \psi_t^{\pi^\theta}(x; I_t(\theta)),
\end{align*}
where $\psi_t^{\pi^\theta}(x; I_t(\theta))$ is given by
\begin{align*}
    &x^\top (Q_t + A_t^\top P_{t+1} A_t) x + 2\left(P_{t+1} W_{t\mid t}^\theta + \sum_{\tau = t+1}^{T-1} \Phi_{\tau, t+1}^\top P_{\tau+1} W_{\tau\mid t}^\theta\right)^\top A_t x + \bar{\psi}_{t+1}(I_t(\theta))\\
    &+ \left(K_t x - \bar{u}_t^\theta(I_t(\theta))\right)^\top (R_t + B_t^\top P_{t+1} B_t) \left(K_t x - \bar{u}_t^\theta(I_t(\theta))\right).
\end{align*}
Using the expected Q function, we know that the optimal policy will pick the action
\begin{align*}
    \pi_t(x; I_t(\theta)) = \argmin_u \mathbb{E}\left[Q_{t}^{\pi^\theta}(x, u; \Xi)\mid I_t(\theta)\right] = - K_t x + \bar{u}_t^\theta(I_t(\theta)).
\end{align*}
Therefore, we see the optimal cost-to-go function at time step $t$ is given by
\begin{align*}
    &\mathbb{E}\left[C_{t}^{\pi^\theta}(x; \Xi)\mid I_t(\theta)\right]\\
    ={}& x^\top Q_t x + (K_t x - \bar{u}_t^\theta(I_t(\theta)))^\top R_t (K_t x - \bar{u}_t^\theta(I_t(\theta)))\\
    &+ ((A_t - B_t K_t)x + B_t \bar{u}_t^\theta(I_t(\theta)))^\top P_{t+1} ((A_t - B_t K_t)x + B_t \bar{u}_t^\theta(I_t(\theta)))\\
    & + 2 \left(P_{t+1} W_{t\mid t}^\theta + \sum_{\tau = t+1}^{T-1} \Phi_{\tau, t+1}^\top P_{\tau + 1} W_{\tau\mid t}^\theta\right)^\top ((A_t - B_t K_t)x + B_t \bar{u}_t^\theta(I_t(\theta))) + \bar{\psi}_{t+1}(I_t(\theta))\\
    ={}& x^\top(Q_t + K_t^\top R_t K_t + (A_t - B_t K_t)^\top P_{t+1} (A_t - B_t K_t)) x - 2\bar{u}_t^\theta(I_t(\theta))^\top R_t K_t x\\
    &+ 2 \bar{u}_t^\theta(I_t(\theta))^\top B_t^\top P_{t+1} (A_t - B_t K_t) x\\
    &+ 2 \left(P_{t+1} W_{t\mid t}^\theta + \sum_{\tau = t+1}^{T-1} \Phi_{\tau, t+1}^\top P_{\tau+1} W_{\tau\mid t}^\theta\right)^\top (A_t - B_t K_t)x\\
    &+ \bar{u}_t^\theta(I_t(\theta))^\top (R_t + B_t^\top P_{t+1} B_t) \bar{u}_t^\theta(I_t(\theta))\\
    &+ 2 \left(P_{t+1} W_{t\mid t}^\theta + \sum_{\tau = t+1}^{T-1} \Phi_{\tau, t+1}^\top P_{\tau+1} W_{\tau\mid t}^\theta\right)^\top B_t \bar{u}_t^\theta(I_t(\theta)) + \bar{\psi}_{t+1}(I_t(\theta)).
\end{align*}
Note that the term $-2\bar{u}_t^\theta(I_t(\theta))^\top R_t K_t x$ and the term $+ 2 \bar{u}_t^\theta(I_t(\theta))^\top B_t^\top P_{t+1} (A_t - B_t K_t) x$ cancel out because $R_t K_t = B_t^\top P_{t+1} (A_t - B_t K_t)$. We also note that the matrix in the first quadratic term can be simplified to
\begin{align*}
    &Q_t + K_t^\top R_t K_t + (A_t - B_t K_t)^\top P_{t+1} (A_t - B_t K_t) \\
    ={}& Q_t + K_t^\top B_t^\top P_{t+1} (A_t - B_t K_t) + (A_t - B_t K_t)^\top P_{t+1} (A_t - B_t K_t) \\
    ={}& Q_t + A_t^\top P_{t+1} (A_t - B_t K_t) \\
    ={}& Q_t + A_t^\top P_{t+1} A_t - A_t^\top P_{t+1} B_t K_t \\
    ={}& Q_t + A_t^\top P_{t+1} A_t - A_t^\top P_{t+1} H_t P_{t+1} A_t \\
    ={}& P_t,
\end{align*}
where the last equation follows by the definition of $P_t$ in \eqref{def:LTV-PKH:e1}.

Therefore, we obtain that
\begin{align*}
    &\mathbb{E}\left[C_{t}^{\pi^\theta}(x; \Xi)\mid I_t(\theta)\right]\\
    ={}& x^\top P_t x + 2 \left((A_t^\top - A_t^\top P_{t+1} H_t)(P_{t+1} W_{t\mid t}^\theta + \sum_{\tau = t+1}^{T-1} \Phi_{\tau, t+1}^\top P_{\tau+1} W_{\tau\mid t}^\theta)\right)^\top x\\
    &+ \bar{u}_t^\theta(I_t(\theta))^\top (R_t + B_t^\top P_{t+1} B_t) \bar{u}_t^\theta(I_t(\theta))\\
    &+ 2 \left(P_{t+1} W_{t\mid t}^\theta + \sum_{\tau = t+1}^{T-1} \Phi_{\tau, t+1}^\top P_{\tau+1} W_{\tau\mid t}^\theta\right)^\top B_t \bar{u}_t^\theta(I_t(\theta)) + \bar{\psi}_{t+1}(I_t(\theta))\\
    ={}& x^\top P_t x + 2 \left(\sum_{\tau = t}^{T-1} \Phi_{\tau, t}^\top P_{\tau+1} W_{\tau\mid t}^\theta\right)^\top x + \bar{\psi}_t(I_t(\theta)),
\end{align*}
where the residual term $\bar{\psi}_t(I_t(\theta))$ is given by
\begin{align*}
    \bar{\psi}_t(I_t(\theta)) ={}& \bar{u}_t^\theta(I_t(\theta))^\top (R_t + B_t^\top P_{t+1} B_t) \bar{u}_t^\theta(I_t(\theta))\\
    &+ 2 \left(P_{t+1} W_{t\mid t}^\theta + \sum_{\tau = t+1}^{T-1} \Phi_{\tau, t+1}^\top P_{\tau+1} W_{\tau\mid t}^\theta\right)^\top B_t \bar{u}_t^\theta(I_t(\theta)) + \bar{\psi}_{t+1}(I_t(\theta)).
\end{align*}
Thus, we have shown the statement of Proposition \ref{thm:LQR-closed-form} and \ref{thm:LQR_closed_form:e1} by induction.
\subsection{Proof of Theorem \ref{thm:pred-power-LTV}}\label{appendix:thm:pred-power-LTV}
By Proposition \ref{thm:LQR-closed-form}, we see that
\begin{align}\label{thm:pred-power-LTV:e1}
    &\mathbb{E}\left[Q_t^{\pi^\theta}(x, u; \Xi) - C_t^{\pi^\theta}(x; \Xi)\mid I_t(\theta) = \iota_t(\theta)\right]\nonumber\\
    ={}& (u - \pi_t^\theta(x; \iota_t(\theta)))^\top (R_t + B_t^\top P_{t+1} B_t) (u - \pi_t^\theta(x; \iota_t(\theta))).
\end{align}
Substituting $u = \bar{\pi}_t(x; \iota_t(\mathbf{0}))$ into the above equation gives that
\begin{subequations}\label{thm:pred-power-LTV:e2}
\begin{align}
    &\mathbb{E}\left[Q_t^{\pi^\theta}(x, \bar{\pi}_t(x; \iota_t(\mathbf{0})); \Xi) - C_t^{\pi^\theta}(x; \Xi)\mid I_t(\theta) = \iota_t(\theta)\right]\nonumber\\
    ={}& (\bar{\pi}_t(x; \iota_t(\mathbf{0})) - \pi_t^\theta(x; \iota_t(\theta)))^\top (R_t + B_t^\top P_{t+1} B_t) (\bar{\pi}_t(x; \iota_t(\mathbf{0})) - \pi_t^\theta(x; \iota_t(\theta)))\nonumber\\
    ={}& (\bar{u}_t^\theta(\iota_t(\theta)) - \bar{u}_t^\zero(\iota_t(\zero)))^\top (R_t + B_t^\top P_{t+1} B_t) (\bar{u}_t^\theta(\iota_t(\theta)) - \bar{u}_t^\zero(\iota_t(\zero))) \label{thm:pred-power-LTV:e2:s1}\\
    ={}& \Tr{(R_t + B_t^\top P_{t+1} B_t) (\bar{u}_t^\theta(\iota_t(\theta)) - \bar{u}_t^\zero(\iota_t(\zero))) (\bar{u}_t^\theta(\iota_t(\theta)) - \bar{u}_t^\zero(\iota_t(\zero)))^\top},\label{thm:pred-power-LTV:e2:s2}
\end{align}
\end{subequations}
where we use the expression of optimal policies in \Cref{thm:LQR-closed-form} in \eqref{thm:pred-power-LTV:e2:s1} and rearrange the terms in \eqref{thm:pred-power-LTV:e2:s2}. Note that by \Cref{thm:LQR-closed-form}, we have
\begin{align*}
    \bar{u}_t^\zero(\iota_t(\zero)) = \mathbb{E}\left[\bar{u}_t^\theta(I_t(\theta))\mid I_t(\zero) = \iota_t(\zero)\right].
\end{align*}
Therefore, by the tower rule and the definition of conditional covariance, we obtain that
\begin{align}\label{thm:pred-power-LTV:e3}
    &\mathbb{E}\left[Q_t^{\pi^\theta}(x, \bar{\pi}_t(x; \iota_t(\mathbf{0})); \Xi) - C_t^{\pi^\theta}(x; \Xi)\mid I_t(\zero) = \iota_t(\zero)\right]\nonumber\\
    ={}& \Tr{(R_t + B_t^\top P_{t+1} B_t) \Cov{\bar{u}_t^\theta(I_t(\theta))\mid I_t(\zero) = \iota_t(\zero)}}.
\end{align}
Let $\{(\bar{X}_t, \bar{U}_t)\}$ denote the (random) trajectory achieved $\bar{\pi}_{0:T-1}$ under problem instance $\Xi$. Since $\bar{X}_t$ is $\mathcal{F}_t(\zero)$-measurable, by \eqref{thm:pred-power-LTV:e3}, we obtain that
\begin{align}\label{thm:pred-power-LTV:e4}
    \mathbb{E}\left[Q_t^{\pi^\theta}(\bar{X}_t, \bar{U}_t; \Xi) - C_t^{\pi^\theta}(\bar{X}_t; \Xi)\mid \mathcal{F}_t(\zero)\right] = \Tr{(R_t + B_t^\top P_{t+1} B_t) \Cov{\bar{u}_t^\theta(I_t(\theta))\mid \mathcal{F}_t(\zero)}},
\end{align}
where we use $\bar{U}_t = \bar{\pi}_t(\bar{X}_t; I_t(\mathbf{0}))$. Note that we have
\begin{align*}
    Q_t^{\pi^\theta}(\bar{X}_t, \bar{U}_t; \Xi) = h_t(\bar{X}_t, \bar{U}_t) + C_{t+1}^{\pi^\theta}(\bar{X}_{t+1}; \Xi).
\end{align*}
Substituting this into \eqref{thm:pred-power-LTV:e4} and taking expectation give that
\begin{align}\label{thm:pred-power-LTV:e5}
    &\mathbb{E}\left[h_t(\bar{X}_t, \bar{U}_t) + C_{t+1}^{\pi^\theta}(\bar{X}_{t+1}; \Xi) - C_t^{\pi^\theta}(\bar{X}_t; \Xi)\right]\nonumber\\
    ={}& \Tr{(R_t + B_t^\top P_{t+1} B_t) \mathbb{E}\left[\Cov{\bar{u}_t^\theta(I_t(\theta))\mid \mathcal{F}_t(\zero)}\right]}.
\end{align}
Summing \eqref{thm:pred-power-LTV:e5} over $t = 0, 1, \ldots, T-1$, we obtain that
\begin{align*}
    \mathbb{E}\left[\sum_{t=0}^{T-1} h_t(\bar{X}_t, \bar{U}_t) - C_0^{\pi^\theta}(\bar{X}_0; \Xi)\right] = \sum_{t=0}^{T-1} \Tr{(R_t + B_t^\top P_{t+1} B_t) \mathbb{E}\left[\Cov{\bar{u}_t^\theta(I_t(\theta))\mid \mathcal{F}_t(\zero)}\right]}.
\end{align*}
Note that the left-hand side equals $P(\theta)$. Thus, we have finished the proof of \Cref{thm:pred-power-LTV}.
\subsection{Proof of the MPC form}\label{appendix:LPC-form}\label{appendix:MPC-form}
In the LQR setting, we can further simplify the MPC policy \eqref{equ:MPC-in-expectation} to be \textit{planning according to $w_{\tau\mid t}^\theta$}:
\begin{equation}\label{equ:MPC-expected-disturbances}
\begin{aligned}
    \argmin_{u_{t:T-1}}\quad
    & \sum_{\tau = t}^{T-1} h_\tau(x_\tau, u_\tau) + h_T(x_T) \\
    \text{s.t.}\quad
    & x_{\tau+1} = f_\tau(x_\tau, u_\tau; w_{\tau\mid t}^\theta),\text{ for }\tau \geq t, \text{ and } x_t = x.
\end{aligned}
\end{equation}

In this section, we show that the MPC policies defined in \eqref{equ:MPC-in-expectation} and \eqref{equ:MPC-expected-disturbances} are equivalent to the optimal policy in Proposition \ref{thm:LQR-closed-form}.

To simplify the notation, we define the large vectors
\begin{align*}
    \vec{x} \coloneqq \begin{bmatrix}
        x_t\\
        x_{t+1}\\
        \vdots\\
        x_T
    \end{bmatrix}, \ \vec{u} \coloneqq \begin{bmatrix}
        u_t\\
        u_{t+1}\\
        \vdots\\
        u_{T-1}
    \end{bmatrix}, \ \text{ and }\vec{w} \coloneqq \begin{bmatrix}
        w_t\\
        w_{t+1}\\
        \vdots\\
        w_{T-1}
    \end{bmatrix}.
\end{align*}
Follow the approach of system level thesis, we know the constraints that
\begin{align*}
    x_{\tau+1} \coloneqq A_\tau x_\tau + B_\tau u_\tau + w_\tau,\text{ for }\tau \geq t, \text{ and } x_t = x
\end{align*}
can be expressed equivalently by the affine relationship
\begin{align*}
    \vec{x} \coloneqq \Phi_x x + \Phi_u \vec{u} + \Phi_w \vec{w}.
\end{align*}
Let $\vec{Q} = \text{Diag}(Q_t, \ldots, Q_{T-1}, P_T)$ and $\vec{R} = \text{Diag}(R_t, \ldots, R_{T-1})$. We know the objective function (with equality constraints)
\begin{align}\label{MPC_form:e1}
    & \sum_{\tau = t}^{T-1} h_\tau(x_\tau, u_\tau) + h_T(x_T)\nonumber\\
    \text{s.t.}\quad
    & x_{\tau+1} = f_\tau(x_\tau, u_\tau; w_t),\text{ for }\tau \geq t, \text{ and } x_t = x,
\end{align}
can be written equivalently in the unconstrained form
\begin{align}\label{MPC_form:e2}
    (\Phi_x x + \Phi_u \vec{u} + \Phi_w \vec{w})^\top \vec{Q} (\Phi_x x + \Phi_u \vec{u} + \Phi_w \vec{w}) + \vec{u}^\top \vec{R} \vec{u}.
\end{align}

We introduce the notations
\begin{align*}
    \vec{W} \coloneqq \begin{bmatrix}
        W_t\\
        W_{t+1}\\
        \vdots\\
        W_{T-1}
    \end{bmatrix}, \ \vec{W}_{\cdot \mid t}^\theta \coloneqq \begin{bmatrix}
        W_{t\mid t}^\theta\\
        W_{t+1\mid t}^\theta\\
        \vdots\\
        W_{T-1\mid t}^\theta
    \end{bmatrix}, \text{ and }\vec{w}_{\cdot \mid t}^\theta \coloneqq \begin{bmatrix}
        w_{t\mid t}^\theta\\
        w_{t+1\mid t}^\theta\\
        \vdots\\
        w_{T-1\mid t}^\theta
    \end{bmatrix}.
\end{align*}
The MPC policy in \eqref{equ:MPC-in-expectation} can be expressed as
\begin{align*}
    \min_{\vec{u}} \Condexp{(\Phi_x x + \Phi_u \vec{u} + \Phi_w \vec{W})^\top \vec{Q} (\Phi_x x + \Phi_u \vec{u} + \Phi_w \vec{W}) + \vec{u}^\top \vec{R} \vec{u}}{I_t(\theta) = \iota_t(\theta)}.
\end{align*}
Because the objective function can be reduced to
\begin{align*}
    &\Condexp{(\Phi_x x + \Phi_u \vec{u} + \Phi_w \vec{W})^\top \vec{Q} (\Phi_x x + \Phi_u \vec{u} + \Phi_w \vec{W}) + \vec{u}^\top \vec{R} \vec{u}}{I_t(\theta) = \iota_t(\theta)}\\
    ={}& (\Phi_x x + \Phi_u \vec{u} + \Phi_w \vec{w}_{\cdot \mid t}^\theta)^\top \vec{Q} (\Phi_x x + \Phi_u \vec{u} + \Phi_w \vec{w}_{\cdot \mid t}^\theta) + \vec{u}^\top \vec{R} \vec{u}\\
    &+ \Condexp{(\Phi_w (\vec{W} - \vec{W}_{\cdot \mid t}^\theta))^\top \vec{Q} \Phi_w (\vec{W} - \vec{W}_{\cdot \mid t}^\theta)}{I_t(\theta) = \iota_t(\theta)},
\end{align*}
where the last term is independent with $x$ and $\vec{u}$. Thus, the MPC policy in \eqref{equ:MPC-in-expectation} is equivalent to
\begin{align*}
    \Condexp{(\Phi_x x + \Phi_u \vec{u} + \Phi_w \vec{W})^\top \vec{Q} (\Phi_x x + \Phi_u \vec{u} + \Phi_w \vec{W}) + \vec{u}^\top \vec{R} \vec{u}}{I_t(\theta) = \iota_t(\theta)},
\end{align*}
which is the MPC policy in \eqref{equ:MPC-expected-disturbances}.

Now, we show that \eqref{equ:MPC-expected-disturbances} is equivalent to the optimal policy in Proposition \ref{thm:LQR-closed-form}. For any sequence $w_{t:T-1}$, let $\texttt{MPC}(x, w_{t:T-1})$ denote the first entry of the solution to
\begin{align}\label{MPC_form:e3}
    \argmin_{u_{t:T-1}}& \sum_{\tau = t}^{T-1} h_\tau(x_\tau, u_\tau) + h_T(x_T)\nonumber\\
    \text{s.t.}\quad
    & x_{\tau+1} = f_\tau(x_\tau, u_\tau; w_t),\text{ for }\tau \geq t, \text{ and } x_t = x,
\end{align}
To show that \eqref{equ:MPC-expected-disturbances} is equivalent to the optimal policy in Proposition \ref{thm:LQR-closed-form}, we only need to show that
\begin{align}\label{MPC_form:e4}
    \texttt{MPC}(x, w_{t:T-1}) = - K_t x - (R_t + B_t^\top P_{t+1} B_t)^{-1} B_t^\top \sum_{\tau = t}^{T-1} \Phi_{\tau+1, t+1}^\top P_{\tau + 1} w_t
\end{align}
holds for any sequence $w_{t:T-1}$. To see this, we consider the case when $w_{t:T-1}$ are deterministic disturbances on and after time step $t$, i.e., the agent knows $w_{t:T-1}$ exactly at time step $t$. In this scenario, we know the optimal policy is to follow the planned trajectory according to MPC in \eqref{MPC_form:e1}. On the other hand, by Proposition \ref{thm:LQR-closed-form}, we know the optimal action to take at time $t$ is $- K_t x - (R_t + B_t^\top P_{t+1} B_t)^{-1} B_t^\top \sum_{\tau = t}^{T-1} \Phi_{\tau+1, t+1}^\top P_{\tau + 1} w_t$. Therefore, the first step planned by MPC must be identical with $- K_t x - (R_t + B_t^\top P_{t+1} B_t)^{-1} B_t^\top \sum_{\tau = t}^{T-1} \Phi_{\tau+1, t+1}^\top P_{\tau + 1} w_t$. Thus, \eqref{MPC_form:e4} holds. And replacing $w_{t:(T-1)}$ with $w_{t:(T-1)\mid t}^\theta$ finishes the proof.
\subsection{Prediction Power Evaluation}\label{appendix:alg:pred-power-eva}
Based on our discussion in \Cref{sec:LQR}, we propose an algorithm (cf. Algorithm \ref{alg:pred-power-eva}) to evaluate the prediction power efficiently given a set of historical problem instances $\{\xi_n\}_{n=1}^N$. Recall that we define the surrogate-optimal action as 
\begin{align}\label{equ:surrogate-optimal-action:appendix}
    \bar{u}_t^*(\Xi) \coloneqq - (R_t + B_t^\top P_{t+1} B_t)^{-1} B_t^\top \sum_{\tau = t}^{T-1} \Phi_{\tau+1, t+1}^\top P_{\tau + 1} W_\tau,
\end{align}
which is the optimal action that an agent should take with the oracle knowledge of all future disturbances at time $t$. In the prediction power given by Theorem \ref{thm:pred-power-LTV}, we can express $\bar{u}_t^\theta(I_t(\theta))$ as $\mathbb{E}\left[\bar{u}_t^*(\Xi)\mid I_t(\theta)\right]$ by \Cref{thm:LQR-closed-form}, which is the expectation of $\bar{u}_t^*(\Xi)$ condition on the the history at time step $t$.

We design \Cref{alg:pred-power-eva} as following: While iterating backward from time step $T-1$ to $0$, the algorithm first constructs a dataset of the surrogate optimal action $\bar{u}_t^*(\Xi)$ as the fitting target. Then, the algorithm estimates the covariance of $\bar{u}_t^*(\Xi)$ when conditioning on $I_t(\mathbf{0})$ and $I_t(\theta)$, respectively, using a subroutine (\Cref{alg:ECCE}). The last step of Algorithm \ref{alg:pred-power-eva} gives the prediction power because $\mathbb{E}\left[\Cov{\bar{u}_t^\theta(I_t(\theta))\mid \mathcal{F}_t(\mathbf{0})}\right]$ can be decomposed as $\mathbb{E}\left[\Cov{\bar{u}_t^{*}(\Xi)\mid I_t(\mathbf{0})}\right] - \mathbb{E}\left[\Cov{\bar{u}_t^{*}(\Xi)\mid I_t(\theta)}\right]$, and we prove this result in \Cref{lemma:cov-cond-expectation}. This decomposition is helpful because otherwise, we would need to evaluate the conditional expectation inside another conditional expectation. Specifically, $\bar{u}_t^\theta(I_t(\theta))$ needs to be approximated by a learned regressor (say, $\phi$) that takes $I_t(\theta)$ as an input. Then, to evaluate $\mathbb{E}\left[\Cov{\bar{u}_t^\theta(I_t(\theta))\mid \mathcal{F}_t(\mathbf{0})}\right]$, we would need to train another regressor to predict the output of $\phi$. Our decomposition avoids this hierarchical dependence.

\begin{algorithm}[t]
    \caption{Prediction Power Evaluation}
    \label{alg:pred-power-eva}
    \begin{algorithmic}[1]
    \REQUIRE Dataset $D$ of problem instances $\{\xi_n\}_{n=1}^N$.\\
    \FOR{$t = T-1, T-2, \ldots, 0$}
    \STATE Compute $P_t, H_t, K_t$ and $\{\Phi_{t, t'}\}_{t' \geq t}$ according to \eqref{def:LTV-PKH:e1} and \eqref{def:LTV-PKH:e1-1}.
    \STATE Compute $M_t = R_t + B_t^\top P_{t+1} B_t$.
    \FOR{$n = 1, 2, \ldots, N$}
        \STATE Compute $\bar{u}_t^*(\xi_n)$ according to \eqref{equ:surrogate-optimal-action} in problem instance $\xi_n.$
    \ENDFOR
    \STATE Call Algorithm \ref{alg:ECCE} to estimate $\Sigma_t^0 \coloneqq \mathbb{E}\left[\Cov{\bar{u}_t^{*}(\Xi)\mid I_t(0)}\right]$ using $\{(\bar{u}_t^{*}(\xi_n), \iota_t^n(0))\}_{n=1}^N$.
    \STATE Call Algorithm \ref{alg:ECCE} to estimate $\Sigma_t^\theta \coloneqq \mathbb{E}\left[\Cov{\bar{u}_t^{*}(\Xi)\mid I_t(\theta)}\right]$ using $\{(\bar{u}_t^{*}(\xi_n), \iota_t^n(\theta))\}_{n=1}^N$.
    \ENDFOR
    \RETURN $P(\theta) = \sum_{t=0}^{T-1} \Tr{\Sigma_t^0 M_t} - \sum_{t=0}^{T-1} \Tr{\Sigma_t^\theta M_t}$
    \end{algorithmic}
\end{algorithm}

\begin{lemma}\label{lemma:cov-cond-expectation}
For any random variable $X$ and two $\sigma$-algebras $\mathcal{F} \subseteq \mathcal{F}'$, the following equation holds
\begin{align*}
    \mathbb{E}\left[\Cov{\mathbb{E}\left[X\mid \mathcal{F}'\right]\mid \mathcal{F}}\right] = \mathbb{E}\left[\Cov{X\mid \mathcal{F}}\right] - \mathbb{E}\left[\Cov{X\mid \mathcal{F}'}\right].
\end{align*}
\end{lemma}
\begin{proof}[Proof of Lemma \ref{lemma:cov-cond-expectation}]
By the law of total covariance, we see that
\begin{align*}
    \Cov{X\middle|\mathcal{F}} = \Cov{\mathbb{E}\left[X\mid \mathcal{F}'\right]\mid \mathcal{F}} + \Condexp{\Cov{X\middle|\mathcal{F}'}}{\mathcal{F}}.
\end{align*}
Taking expectation on both sides gives that
\begin{align*}
    \mathbb{E}\left[\Cov{X\mid \mathcal{F}}\right] = \mathbb{E}\left[\Cov{\mathbb{E}\left[X\mid \mathcal{F}'\right]\mid \mathcal{F}}\right] + \mathbb{E}\left[\Cov{X\mid \mathcal{F}'}\right],
\end{align*}
which is equivalent to the statement of Lemma \ref{lemma:cov-cond-expectation}.
\end{proof}

\textbf{Evaluation of the Expected Conditional Covariance.} For two general random variables $X$ and $Y$, we follow a standard procedure to evaluate the expectation of their conditional covariance $\mathbb{E}\left[\Cov{Y\mid X}\right]$ using a dataset $\{(x_n, y_n)\}$ that is independently sampled from the joint distribution of $(X, Y)$ (Algorithm \ref{alg:ECCE}). The algorithm first train a regressor $\psi$ that approximates the conditional expectation $\mathbb{E}\left[X\mid Y\right]$, where we use the definition:
\begin{align*}
    \mathbb{E}\left[Y\mid X\right] = \min_{\psi \text{ is any function.}} \mathbb{E}\left[\norm{Y - \psi(X)}_2^2\right].
\end{align*}
Then, $\psi$ is used for evaluating the conditional covariance. During training, we split the dataset to the train, validation, and test datasets in order to prevent overfitting.

\begin{algorithm}
    \caption{Expected Conditional Covariance Estimator (ECCE)}\label{alg:ECCE}
    \begin{algorithmic}[1]
    \REQUIRE Dataset $D$ that consists input/output pair $(x_n, y_n)$.
    \STATE Split the dataset $D$ to $D_\text{train}, D_\text{val},$ and $D_\text{test}$.
    \STATE Initialize a regressor $\psi$ with input $x$ and target output $y$.
    \STATE Fit $\psi$ to $D_\text{train}$ with MSE and use $D_\text{val}$ to prevent over-fit.
    \RETURN $\Sigma \coloneqq \frac{1}{\abs{D_\text{test}}}\sum_{n \in D_\text{test}} (y_n - \psi(x_n))(y_n - \psi(x_n))^\top$.
    \end{algorithmic}
\end{algorithm}
\subsection{Details of Examples in Section \ref{sec:main:power-vs-accuracy}}\label{appendix:example-details}
In this section, we present the specific instantiation of \Cref{example:one-step-dim-rotation} in \Cref{sec:main:power-vs-accuracy} and another example (\Cref{example:multi-step-1d}) for the mismatch between prediction power and prediction accuracy. The code for the experiments can be found at \url{https://github.com/yihenglin97/Prediction-Power}.
\subsubsection{Instantiation of Example \ref{example:one-step-dim-rotation}}\label{appendix:example-details:double-integrator}
\begin{figure*}
\centering
\begin{minipage}{.48\textwidth}
  \centering
  \includegraphics[width=.95\linewidth]{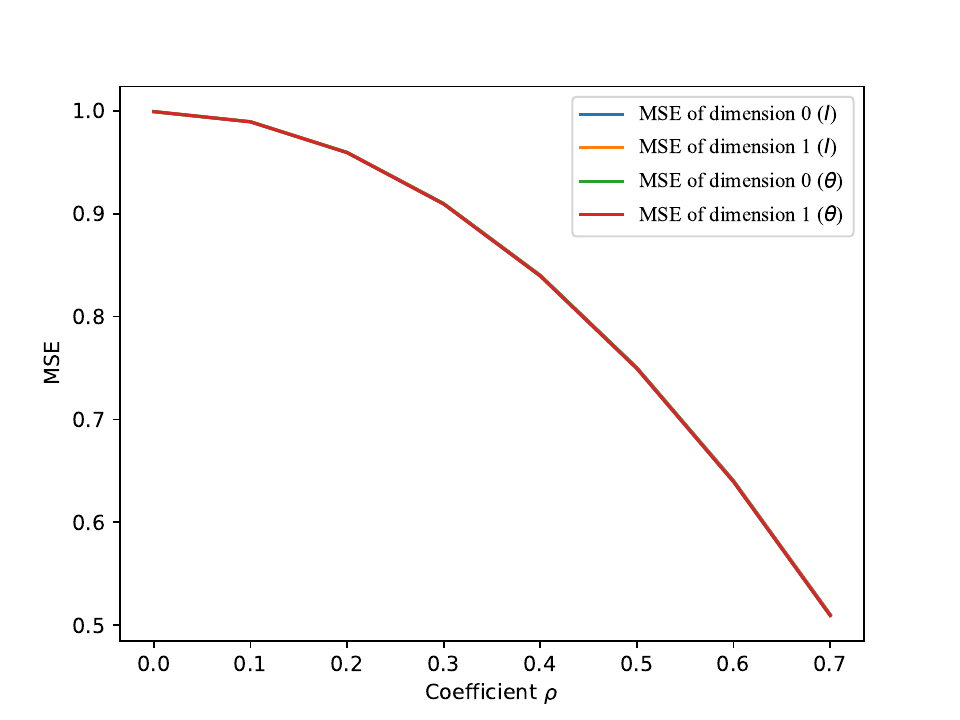}
  \caption{Example \ref{example:one-step-dim-rotation}: MSE - $\rho$ curve.}
  \label{fig:1_step_mse}
\end{minipage}%
\hfill
\begin{minipage}{.48\textwidth}
  \centering
  \includegraphics[width=.95\linewidth]{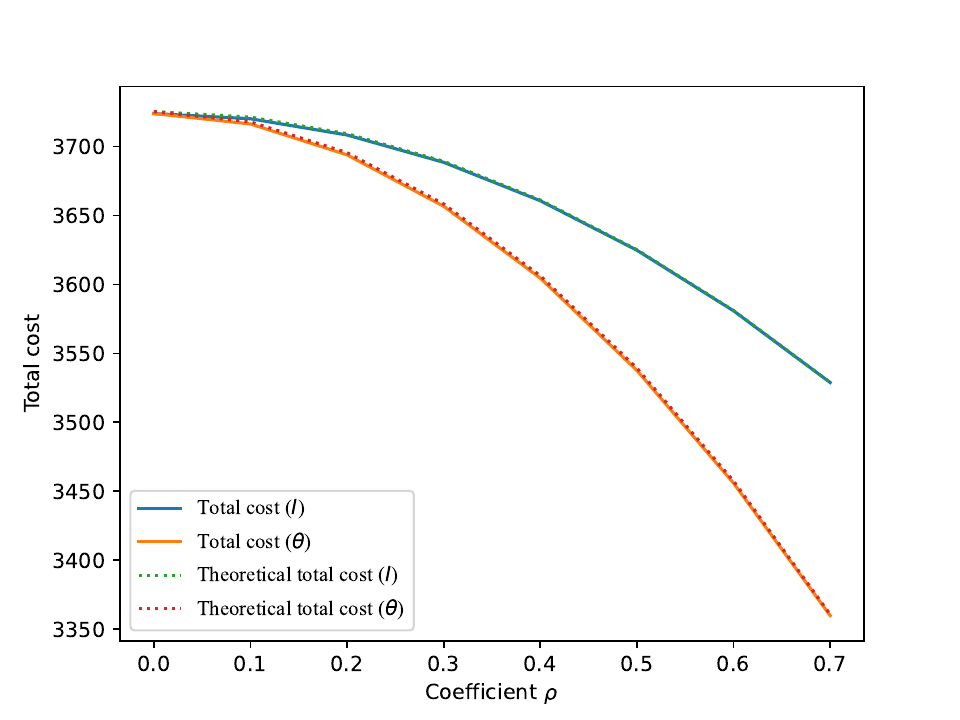}
  \caption{Example \ref{example:one-step-dim-rotation}: Control cost - $\rho$ curve.}
  \label{fig:1_step_cost_comparison}
\end{minipage}
\end{figure*}

We instantiate Example \ref{example:one-step-dim-rotation} with the following parameters:
\begin{align*}
    A = \begin{bmatrix}
        1 & 0.1\\
        0 & 1
    \end{bmatrix}, \ B = \begin{pmatrix}
        0\\
        0.1
    \end{pmatrix}, \ Q = \begin{pmatrix}
        1 & \\
         & 1
    \end{pmatrix}, \ R = (1), \text{ and }\theta \coloneqq \begin{bmatrix}
        1 & 0.99\\
        0 & 0.141
    \end{bmatrix}.
\end{align*}
Under different values of coefficient $\rho$, we train a linear regressor to predict each entry of $W_t$ from $V_t(\theta)$ (or $V_t(I)$) over a train dataset with $64000$ independent samples. We plot in the MSE - $\rho$ curve on a test dataset with $16000$ independent samples in Figure \ref{fig:1_step_mse}. From the plot, we see that the predictors $V_t(\theta)$ and $V_t(I)$ achieve the same MSE when predicting each entry of $W_t$ under each $\rho \in \{0, 0.1, \ldots, 0.7\}$.

Then, we use the trained linear regressors as $W_{t\mid t}^\theta$ and $W_{t\mid t}^I$ to implement the optimal policy in Proposition \ref{thm:LQR-closed-form}. We plot the averaged total cost over $16000$ trajectories with horizon $T = 100$ in Figure \ref{fig:1_step_cost_comparison}. From the plot, we see that the optimal policies under the predictors $V_t(\theta)$ and $V_t(I)$ achieve significantly different control costs when $\rho > 0$. We also plot the theoretical expected control cost in Figure \ref{fig:1_step_cost_comparison} to verify this cost difference. Running this experiment takes about 50 seconds on Apple Mac mini with Apple M1 CPU.

\subsubsection{An One-dimension Example}
\begin{figure}
    \centering
    \includegraphics[width=0.5\linewidth]{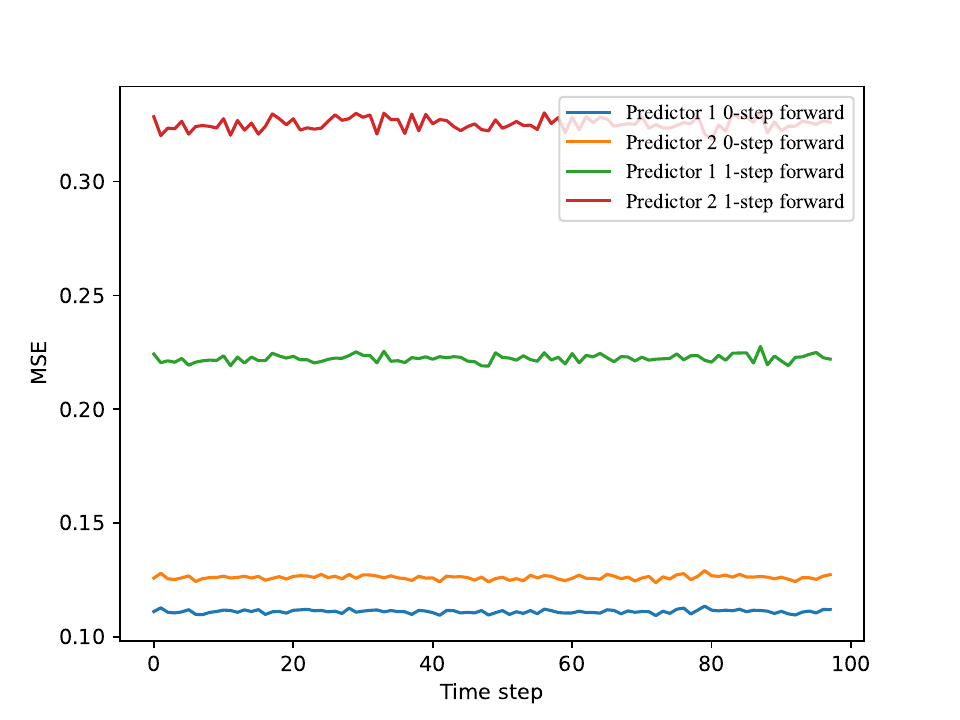}
    \caption{Example \ref{example:multi-step-1d}: MSE - time curve.}
    \label{fig:multi_step_mse}
\end{figure}

We also provide an example with $n = 1$, where the prediction $V_t(\theta)$ is correlated with two steps of future disturbances $W_t$ and $W_{t+1}$.
\begin{example}\label{example:multi-step-1d}
Suppose the disturbance at each time step can be decomposed as $W_t = \sum_{i=0}^2 W_t^{(i)}$, where the $\{W_t^{(i)}\}_{i=0}^2$ are independently sampled from three mean-zero distributions. We compare two predictors: $V_t(1) = \left(W_t^{(1)}, W_{t+1}^{(0)}\right)$ and $V_t(2) = P \left(W_t^{(0)} + W_t^{(1)}\right) + (A^\top - A^\top P H) P W_{t+1}^{(0)}$. They have the same prediction power when used in the control problem because
\begin{align*}
    \bar{u}_t^2(I_t(2)) = P \left(W_t^{(0)} + W_t^{(1)}\right) + (A^\top - A^\top P H) P W_{t+1}^{(1)} = \bar{u}_t^1(I_t(1)).
\end{align*}
However, we know that $\mathcal{F}_t(1)$ is a strict super set of $\mathcal{F}_t(2)$, thus $V_t(1)$ can achieve a better MSE than $V_t(2)$ when predicting the disturbances. This is empirically verified in a 1D LQR problem with $A = B = Q = R = (1)$ and $W_t^{(i)} \overset{\text{i.i.d.}}{\sim} N(0, 1)$, as we plot in Figures \ref{fig:multi_step_mse}. In the simulation, we train linear regressors to predict $W_t$ and $W_{t+1}$ with the history $I_t(1)$ or $I_t(2)$ for each time step $t < T = 100$ over a train dataset of size $160000$. Then, we plot the MSE - time curve on a test dataset of size $40000$. Running this experiment takes about 270 seconds on Apple Mac mini with Apple M1 CPU.
\end{example}

\subsection{Details of Example \ref{example:pred-power-and-policy-opt}}\label{appendix:example-pred-power-and-policy-opt}
We instantiate \Cref{example:pred-power-and-policy-opt} with the same dynamics and costs as \Cref{example:one-step-dim-rotation}, i.e.,
\begin{align*}
    A = \begin{bmatrix}
        1 & 0.1\\
        0 & 1
    \end{bmatrix}, \ B = \begin{pmatrix}
        0\\
        0.1
    \end{pmatrix}, \ Q = \begin{pmatrix}
        1 & \\
         & 1
    \end{pmatrix}, \text{ and } R = (1).
\end{align*}
To build the predictors, we sample the true disturbance $W_t \overset{\text{i.i.d.}}{\sim} N(0, I)$ and fix the coefficient $\rho = 0.5$. The online policy optimization starts with the initial policy parameter $\Upsilon_0 = \mathbf{0}$. When implementing M-GAPS in both scenarios, we use the decaying learning rate sequence $\eta_t = (1 + t/1000)^{-0.5}$. The optimal predictive policy for using $V_{0:t-1}(1)$ or $V_{0:t-1}(2)$ are $\pi^1_{0:T-1}$ and $\pi^2_{0:T-1}$, whose closed-form expressions are given by \Cref{thm:LQR-closed-form}. Note that for the history $\iota_t(1)$, the optimal predictive policy $\pi_t^1$ only depends on $v_t(1)$ because all other entries are independent with future disturbances $W_{t:T-1}$. Similarly, for the history $\iota_t(2)$, the optimal predictive policy $\pi_t^2$ only depends on $v_{t-1}(2)$ and $v_t(2)$.

In Figures \ref{fig:Pred-power-vs-avg-cost:Current} and \ref{fig:Pred-power-vs-avg-cost:Next}, we compute the average cost improvement of M-GAPS (or the optimal predictive policy) against the optimal no-prediction controller $\bar{\pi}_t(x) = - K x$. That is, on each problem instance $\xi$, we plot
\begin{align*}
    \frac{1}{t+1}\left(-(\text{cost of M-GAPS until time }t) + (\text{cost of }\bar{\pi}\text{ until time }t)\right)
\end{align*}
for time $t = 0, 1, \ldots, T-1$. The prediction power (averaged over time) is given by $P(\theta)/T$. We simulate 30 random trajectories with $T = 80000$ and plot the mean with the 25-th and 75-th percentiles as shaded areas. From the plots, we see that M-GAPS' average cost improvement converges towards the prediction power over time in the first scenario but stays far away with the prediction power in the second scenario. This is as expected, because the optimal predictive policy $\pi_t^2$ is not in the candidate policy set of M-GAPS in the second scenario. Simulating the first scenario takes about 200 seconds on Apple Mac mini with Apple M1 CPU. The second takes about 210 seconds on the same hardware.

\section{Proof of Theorem \ref{thm:Bellman-improvement-not-worse}}\label{appendix:thm:Bellman-improvement-not-worse}
Since we assume $x_0$ is the initial state (deterministic) and $\pi^\theta$ is the optimal policy under the predictor with parameter $\theta$, we have
\begin{align*}
    \bE{C_0^{\pi^\theta}(x_0; \Xi)} = J^{\pi^\theta}(\theta) = J^*(\theta).
\end{align*}
Similarly, we also have that
\begin{align*}
    \bE{C_0^{\bar{\pi}}(x_0; \Xi)} = J^{\bar{\pi}}(\mathbf{0}) = J^*(\mathbf{0}).
\end{align*}

Let $\{\bar{X}_{0:T}, \bar{U}_{0:T-1}\}$ be the trajectory of the baseline controller $\bar{\pi}_{0:T-1}$ under instance $\Xi$ starting from $\bar{X}_0 = x_0$. First, we will prove by backwards induction that the difference in cumulative costs between the optimal controller $\pi^\theta$ and $\bar\pi$ has the following decomposition:
\begin{equation}\label{thm:monotonic-improvement:e1}
    C_0^{\pi^\theta}(x_0; \Xi) - C_0^{\bar{\pi}}(x_0; \Xi) = \sum_{t = 0}^{T-1} \left(C_t^{\pi^\theta}(\bar{X}_t; \Xi) - Q_t^{\pi^\theta}(\bar{X}_t, \bar{U}_t; \Xi)\right).
\end{equation}
For the base case at time $T-1$, we apply the definition of $C_{T-1}^{\bar\pi}$ to get
\[
    C_{T-1}^{\pi^\theta}(\bar{X}_{T-1}; \Xi) - C_{T-1}^{\bar\pi}(\bar{X}_{T-1}; \Xi)
    = C_{T-1}^{\pi^\theta}(\bar{X}_{T-1}; \Xi) - Q_{T-1}^{\pi^\theta}(\bar{X}_{T-1}, \bar{U}_{T-1}; \Xi).
\]
For the inductive step, suppose that
\[
    C_{\tau+1}^{\pi^\theta}(\bar{X}_{\tau+1}; \Xi) - C_{\tau+1}^{\bar\pi}(\bar{X}_{\tau+1}; \Xi) = \sum_{t=\tau+1}^{T-1} \left(C_t^{\pi^\theta}(\bar{X}_t; \Xi) - Q_t^{\pi^\theta}(\bar{X}_t, \bar{U}_t; \Xi)\right).
\]
Note that for any $t < T$,
\[
    Q_t^{\bar\pi}(\bar{X}_t, \bar{U}_t; \Xi)
    = Q_t^{\pi^\theta}(\bar{X}_t, \bar{U}_t; \Xi) - \left( C_{t+1}^{\pi^\theta}(\bar{X}_{t+1}; \Xi) - C_{t+1}^{\bar\pi}(\bar{X}_{t+1}; \Xi) \right).
\]
Therefore,
\begin{align*}
    & C_\tau^{\pi^\theta}(\bar{X}_\tau; \Xi) - C_\tau^{\bar\pi}(\bar{X}_\tau; \Xi) \\
    &= C_\tau^{\pi^\theta}(\bar{X}_\tau; \Xi) - Q_\tau^{\bar\pi}(\bar{X}_\tau, \bar{U}_\tau; \Xi) \\
    &= C_\tau^{\pi^\theta}(\bar{X}_\tau; \Xi) - \left[ Q_\tau^{\pi^\theta}(\bar{X}_\tau, \bar{U}_\tau; \Xi) - \left( C_{\tau+1}^{\pi^\theta}(\bar{X}_{\tau+1}; \Xi) - C_{\tau+1}^{\bar\pi}(\bar{X}_{\tau+1}; \Xi) \right) \right] \\
    &= C_\tau^{\pi^\theta}(\bar{X}_\tau; \Xi) - Q_\tau^{\pi^\theta}(\bar{X}_\tau, \bar{U}_\tau; \Xi) + \sum_{t=\tau+1}^{T-1} \left(C_t^{\pi^\theta}(\bar{X}_t; \Xi) - Q_t^{\pi^\theta}(\bar{X}_t, \bar{U}_t; \Xi)\right) \\
    &= \sum_{t=\tau}^{T-1} \left(C_t^{\pi^\theta}(\bar{X}_t; \Xi) - Q_t^{\pi^\theta}(\bar{X}_t, \bar{U}_t; \Xi)\right).
\end{align*}
This completes the induction.

Next, define $U_t \coloneq \pi_t^\theta(\bar{X}_t; I_t(\theta))$. Note that $U_t$ is $\mathcal{F}_t(\theta)$-measurable, and $\bar{U}_t$ is $\mathcal{F}_t(\zero)$-measurable and therefore also $\mathcal{F}_t(\theta)$-measurable. Because we assume the matrices $M_{0:T-1}$ satisfy \Cref{cond:Q-function-quadratic-growth},
\begin{equation}\label{thm:monotonic-improvement:e5}
    \mathbb{E}\left[C_t^{\pi^\theta}(\bar{X}_t; \Xi) \mid I_t(\theta)\right] \leq \mathbb{E}\left[Q_t^{\pi^\theta}(\bar{X}_t, \bar{U}_t; \Xi) \mid I_t(\theta)\right] - \Tr{M_t (\bar{U}_t - U_t)(\bar{U}_t - U_t)^\top}.
\end{equation}

Let $\tilde{U}_t \coloneqq \mathbb{E}\left[U_t\mid I_t(\zero)\right]$. We see that
\begin{subequations}\label{thm:monotonic-improvement:e5-1}
\begin{align}
    &\mathbb{E}\left[(\bar{U}_t - U_t)(\bar{U}_t - U_t)^\top\mid I_t(\zero)\right]\nonumber\\
    ={}& \mathbb{E}\left[(\tilde{U}_t - U_t)(\tilde{U}_t - U_t)^\top\mid I_t(\zero)\right] + \mathbb{E}\left[(\tilde{U}_t - U_t)(\bar{U}_t - \tilde{U}_t)^\top\mid I_t(\zero)\right]\nonumber\\
    &+ \mathbb{E}\left[(\bar{U}_t - \tilde{U}_t)(\tilde{U}_t - U_t)^\top\mid I_t(\zero)\right] + \mathbb{E}\left[(\bar{U}_t - \tilde{U}_t)(\bar{U}_t - \tilde{U}_t)^\top\mid I_t(\zero)\right]\nonumber\\
    ={}& \Cov{\pi_t^\theta(\bar{X}_t; I_t(\theta))\mid I_t(\zero)} + \mathbb{E}\left[\tilde{U}_t - U_t\mid I_t(\zero)\right](\bar{U}_t - \tilde{U}_t)^\top\nonumber\\
    &+ (\bar{U}_t - \tilde{U}_t)\mathbb{E}\left[\tilde{U}_t - U_t\mid I_t(\zero)\right]^\top + (\bar{U}_t - \tilde{U}_t)(\bar{U}_t - \tilde{U}_t)^\top\label{thm:monotonic-improvement:e5-1:s1}\\
    ={}& \Cov{\pi_t^\theta(\bar{X}_t; I_t(\theta))\mid I_t(\zero)} + (\bar{U}_t - \tilde{U}_t)(\bar{U}_t - \tilde{U}_t)^\top,\label{thm:monotonic-improvement:e5-1:s2}
\end{align}
\end{subequations}
where we use $(\bar{U}_t - \tilde{U}_t)$ is $\mathcal{F}_t(\zero)$-measurable in \eqref{thm:monotonic-improvement:e5-1:s1}; we use the definition of $\tilde{U}_t$ in \eqref{thm:monotonic-improvement:e5-1:s2}.

Applying the towering rule in \eqref{thm:monotonic-improvement:e1} and substituting in \eqref{thm:monotonic-improvement:e5} gives that
\begin{align}\label{thm:monotonic-improvement:e6}
    \mathbb{E}\left[C_0^{\pi^\theta}(x_0; \Xi) - C_0^{\bar{\pi}}(x_0; \Xi)\right] ={}& \sum_{t = 0}^{T-1} \mathbb{E}\left[C_t^{\pi^\theta}(\bar{X}_t; \Xi) - Q_t^{\pi^\theta}(\bar{X}_t, \bar{U}_t; \Xi)\right]\nonumber\\
    ={}& \sum_{t = 0}^{T-1} \mathbb{E}\left[\mathbb{E}\left[C_t^{\pi^\theta}(\bar{X}_t; \Xi)\mid I_t(\theta)\right] - \mathbb{E}\left[Q_t^{\pi^\theta}(\bar{X}_t, \bar{U}_t; \Xi)\mid I_t(\theta)\right]\right]\nonumber\\
    \leq{}& - \sum_{t = 0}^{T-1} \mathbb{E}\left[\Tr{M_t (\bar{U}_t - U_t)(\bar{U}_t - U_t)^\top}\right],\nonumber\\
    ={}& - \sum_{t = 0}^{T-1} \Tr{M_t \mathbb{E}\left[(\bar{U}_t - U_t)(\bar{U}_t - U_t)^\top\right]}.
\end{align}

If the stronger Condition \ref{cond:optimal-policy-psd-covariance}\,\ref{cond:optimal-policy-psd-covariance:a} holds, by \eqref{thm:monotonic-improvement:e5-1}, since $\bar{X}_t$ is $\mathcal{F}_t(\zero)$-measurable, we have
\begin{align}\label{thm:monotonic-improvement:e5-2}
    \mathbb{E}\left[(\bar{U}_t - U_t)(\bar{U}_t - U_t)^\top\right]
    &= \mathbb{E}\left[\mathbb{E}\left[(\bar{U}_t - U_t)(\bar{U}_t - U_t)^\top\mid I_t(\zero)\right]\right] \nonumber\\
    &\succeq \mathbb{E}\left[\Cov{\pi_t^\theta(\bar{X}_t; I_t(\theta))\mid I_t(\zero)}\right] \succeq \Sigma_t.
\end{align}
Then, we can apply \eqref{thm:monotonic-improvement:e5-2} in \eqref{thm:monotonic-improvement:e6} to obtain that
\begin{align}\label{thm:monotonic-improvement:e7}
    \mathbb{E}\left[C_0^{\pi^\theta}(x_0; \Xi) - C_0^{\bar{\pi}}(x_0; \Xi)\right] \leq - \sum_{t=0}^{T-1}\Tr{M_t \Sigma_t}.
\end{align}

Else, if the weaker Condition \ref{cond:optimal-policy-psd-covariance}\,\ref{cond:optimal-policy-psd-covariance:b} holds, by \eqref{thm:monotonic-improvement:e5-1}, since $\bar{X}_t$ is $\mathcal{F}_t(\zero)$-measurable, we have
\begin{align}\label{thm:monotonic-improvement:e5-3}
    \Tr{\mathbb{E}\left[(\bar{U}_t - U_t)(\bar{U}_t - U_t)^\top\right]} ={}& \mathbb{E}\left[\Tr{\mathbb{E}\left[(\bar{U}_t - U_t)(\bar{U}_t - U_t)^\top\mid I_t(\zero)\right]}\right]\nonumber\\
    \geq{}& \mathbb{E}\left[\Tr{\Cov{\pi_t^\theta(\bar{X}_t; I_t(\theta))\mid I_t(\zero)}}\right] \geq \sigma_t.
\end{align}
Note that for any positive semi-definite matrices $A, B, C$ such that $A \succeq C \succeq 0$, we have
\begin{align*}
    \Tr{AB}=\Tr{CB}+\Tr{(A-C)B} \geq \Tr{CB}.
\end{align*}
Since $M_t \succeq \mu_\text{min}(M_t) I$, we can apply \eqref{thm:monotonic-improvement:e5-3} in \eqref{thm:monotonic-improvement:e6} to obtain that
\begin{align*}
    \mathbb{E}\left[C_0^{\pi^\theta}(x_0; \Xi) - C_0^{\bar{\pi}}(x_0; \Xi)\right] \leq{}& -\sum_{t=0}^{T-1}\text{Tr}\left\{\mu_\text{min}(M_t) I\cdot \mathbb{E}\left[(\bar{U}_t - U_t)(\bar{U}_t - U_t)^\top\right]\right\}\\
    \leq{}&  - \sum_{t=0}^{T-1}\mu_\text{min}(M_t) \sigma_t.
\end{align*}

\section{Proofs for LTV Dynamics with General Costs}\label{appendix:LTV-general-costs}
In this section, we first provide a proof outline of \Cref{coro:one-step-pred-power-lower} (\Cref{appendix:LTV-general-costs-outline}). Then, we discuss an example where the MPC in \eqref{equ:MPC-in-expectation} is suboptimal (\Cref{appendix:MPC-counterexample}). Lastly, we provide the proofs for the key technical results required by the proof of \Cref{coro:one-step-pred-power-lower}.
\subsection{Proof Outline of Theorem \ref{coro:one-step-pred-power-lower}}\label{appendix:LTV-general-costs-outline}
Assumption \ref{assump:well-conditioned-cost} makes two requirements about the well-conditioned cost functions, which are standard in the literature of online optimization and control \citep{lin2021perturbation,lin2022bounded}. For the last requirement, we additionally require $\ell_A < 1$, which implies that the system is open-loop stable. Under Assumption \ref{assump:well-conditioned-cost}, the expected cost-to-go function is a well-conditioned function, which is important for establishing Conditions \ref{cond:Q-function-quadratic-growth} and \ref{cond:optimal-policy-psd-covariance}\,\ref{cond:optimal-policy-psd-covariance:b}. 
We state this result formally in Lemma \ref{lemma:well-condition-assump-verification-strongly-convex} in \Cref{sec:LTV-well-conditioned}, which establishes uniform bounds for the strongly convexity/smoothness of the conditional expectation of cost-to-go functions: $\mu_t$ is uniformly bounded below by $\mu_x$ and $\ell_t$ is uniformly bounded above by $\frac{\ell_x}{1 - \ell_A}$.
We present a proof sketch of \Cref{lemma:well-condition-assump-verification-strongly-convex} and defer the formal proof to \Cref{appendix:lemma:well-condition-assump-verification-strongly-convex}.

Starting from time step $T$, we know the cost-to-go $C_T^{\pi^\theta}(x; \Xi)$ equals to the terminal cost $h_t^x(x)$. It satisfies the strong convexity/smoothness directly by Assumption \ref{assump:well-conditioned-cost}. We repeat the following induction iterations: Given $\mathbb{E}\left[C_{t+1}^{\pi^\theta}(x; \Xi)\mid I_{t+1}(\theta)\right]$ at time $t+1$, we define an auxiliary function that adds in the disturbance residual $W_t - W_{t\mid t}^\theta$ and condition on the history at time $t$:
\begin{align}\label{equ:auxiliary-cost-to-go-function}
    \bar{C}_{t+1}^{\pi^\theta}(x; \iota_t(\theta)) \coloneqq \mathbb{E}\left[C_{t+1}^{\pi^\theta}(x + W_t - W_{t\mid t}^\theta; \Xi)\mid I_t(\theta) = \iota_t(\theta)\right].
\end{align}
It can be expressed as $\mathbb{E}\left[\left.\mathbb{E}\left[C_{t+1}^{\pi^\theta}(x + W_t - W_{t\mid t}^\theta; \Xi)\mid I_{t+1}(\theta)\right]\right| I_t(\theta) = \iota_t(\theta)\right]$ by the tower rule. Thus, we know function $\bar{C}_{t+1}^{\pi^\theta}$ is strongly convex and smooth in $x$ because these properties are preserved after taking the expectation. Then, we can obtain the expected cost-to-go function $\mathbb{E}\left[C_t^{\pi^\theta}(x; \Xi)\mid I_t(\theta) = \iota_t(\theta)\right] = h_t^x(x) + \min_u\left(h_t^u(u) + \bar{C}_{t+1}^{\pi^\theta}(A_t x + B_t u + w_{t\mid t}^\theta; \iota_t(\theta))\right).$ We use an existing tool called \textit{infimal convolution} to study the optimal value of the this optimization problem as a function of $x$. Specifically, define an operator $\square_B$:\footnote{If $\omega$ takes an additional parameter $w$, we denote $\left(f \square_B \omega\right)(x; w) \coloneqq \min_{u \in \mathbb{R}^m} \left\{f(u) + \omega(x - B u; w)\right\}$}
\begin{align}\label{equ:infimal-conv-variant-def}
    \left(f \square_B \omega\right)(x) \coloneqq \min_{u \in \mathbb{R}^m} \left\{f(u) + \omega(x - B u)\right\} \text{ for } f: \mathbb{R}^m \to \mathbb{R} \text{ and }\omega: \mathbb{R}^n \to \mathbb{R}.
\end{align}
One can show that if $f$ and $\omega$ are well-conditioned functions, then $\left(f \square_B \omega\right)$ is also well-conditioned (see Appendix \ref{appendix:lemma:infimal-conv-variant} for the formal statement and proof). We can use this result to show the expected cost-to-go function $\mathbb{E}\left[C_t^{\pi^\theta}(x; \Xi)\mid I_t(\theta) = \iota_t(\theta)\right] = h_t^x(x) + (h_t^u\square_{(-B_t)}\bar{C}_{t+1}^{\pi^\theta})(A_t x + w_{t\mid t}^\theta; \iota_t(\theta)),$ is also well-conditioned in $x$ at time step $t$, which completes the induction.

For the second condition on the covariance of $\pi_t^\theta$'s actions, we note that $\lambda_t(\theta)$ in \Cref{assump:one-step-prediction-indep-pairs} should be positive as long as $V_t(\theta)$ has some weak correlation with $W_t$. Under Assumption \ref{assump:one-step-prediction-indep-pairs}, we can express the optimal policy as
\begin{align}\label{equ:joint-gaussian-intuition:e2}
    \pi_t^\theta(x; I_t(\theta)) &\coloneqq \argmin_u \left(h_t^u(u) + \bar{C}_{t+1}^{\pi^\theta}(A_t x + B_t u + W_{t\mid t}^\theta)\right).
\end{align}
While the original definition of $\bar{C}_{t+1}^{\pi^\theta}$ in \eqref{equ:auxiliary-cost-to-go-function} requires the history $\iota_t(\theta)$ as an input, it no longer depends on the history under Assumption \ref{assump:one-step-prediction-indep-pairs}. We defer the proof to Appendix \ref{appendix:coro:one-step-pred-power-lower}.

We can express $\pi_t^\theta(x; I_t(\theta))$ as the solution to $(h_t^u\square_{(-B_t)}\bar{C}_{t+1}^{\pi^\theta})(A_t x + W_{t\mid t})$. For some distributions including Gaussian, the covariance in the input of an infimal convolution will be passed through to its optimal solution. Specifically, let $u_{(f\square_B \omega)}(x)$ denote the solution to the optimization problem \eqref{equ:infimal-conv-variant-def}. When $\omega$ and $f$ are well-conditioned, we can derive a lower bound on the trace of the covariance $\Tr{\Cov{u_{(f\square_B \omega)}(X)}}$ that depends on the covariance of $X$. Due to space limit, we defer the formal statement of this result and its proof to Lemma \ref{lemma:infimal-conv-opt-solution-perturb-lower} in Appendix \ref{appendix:lemma:infimal-conv-variant}. Using this property and the observation that $\pi_t^\theta(x; I_t(\theta))$ can be expressed as $u_{(h_t^u\square_{-B_t}\bar{C}_{t+1}^{\pi^\theta})}(A_t x + W_{t\mid t}^\theta)$, we can directly verify that Condition \ref{cond:optimal-policy-psd-covariance}\,\ref{cond:optimal-policy-psd-covariance:b} holds with
\begin{align}\label{equ:one-step-pred-second-cond:appendix}
    \Tr{\Cov{\pi_t^\theta(x; I_t(\theta))\mid \mathcal{F}_t(0)}} \geq \sigma_t \coloneqq \frac{n \lambda_t(\theta) \mu_{t+1}^2 \cdot \mu_B}{2(\ell_u + \ell_{t+1}\sqrt{\ell_B})^2}.
\end{align}
Since Lemma \ref{lemma:well-condition-assump-verification-strongly-convex} and \eqref{equ:one-step-pred-second-cond:appendix} imply that Conditions \ref{cond:Q-function-quadratic-growth} and \ref{cond:optimal-policy-psd-covariance}\,\ref{cond:optimal-policy-psd-covariance:b} hold with $M_t = \mu_t I$ and $\sigma_t$ respectively, we can apply \Cref{thm:Bellman-improvement-not-worse} to obtain the prediction power lower bound in \Cref{coro:one-step-pred-power-lower}.
\subsection{Example: MPC can be suboptimal}\label{appendix:MPC-counterexample}
We first highlight the challenge by showing that MPC can be suboptimal, i.e., only planning and optimizing based on the current information might be suboptimal when the cost functions are not quadratic.

Consider a 2-step optimal control problem (1-dimension):
\begin{align*}
    X_1 = X_0 + U_0, \text{ and } X_2 = X_1 + U_1 + W_1.
\end{align*}
The cost functions are given by
\begin{align*}
    h_0(x, u) = x^2 + u^2,\ h_1(x, u) = x^2 + u^2, \text{ and }h_2(x) = \begin{cases}
        x^2, & \text{ if } x \leq 0,\\
        +\infty, & \text{ otherwise.}
    \end{cases}
\end{align*}
Suppose $W_1$ is a random variable that satisfies $\mathbb{P}(W_1 = 1) = p$ and $\mathbb{P}(W_1 = 0) = 1 - p$, where $0 < p < 1$. At time $0$, we don't have any knowledge about $W_1$ (i.e., $W_1$ is independent with $I_0(\theta)$). However, at time $1$, we can predict $W_1$ exactly, which means $\sigma(W_1) \subseteq \mathcal{F}_1(\theta)$.

Suppose the system starts at $x_0 = 0$. At time step $0$, MPC \eqref{equ:MPC-in-expectation} solves the optimization
\begin{align}\label{MPC-suboptimal-example:e0}
    \min_{u_0, u_1}& \mathbb{E}\left[h_0(X_0, u_0) + h_1(X_1, u_1) + h_2(X_2)\mid I_0(\theta)\right]\nonumber\\
    \text{ s.t. }&X_0 = 0,\ X_1 = X_0 + u_0, \ X_2 = X_1 + u_1 + W_1.
\end{align}
Since $I_0(\theta)$ is independent with $W_1$, the optimization problem can be expressed equivalently as
\begin{align*}
    &\min_{u_0, u_1} u_0^2 + (u_0^2 + u_1^2) + \mathbb{E}\left[h_2(u_0 + u_1 + W_1)\right]\\
    ={}& \min_{u_0, u_1} 2 u_0^2 + u_1^2 + 1, \text{ s.t. } u_0 + u_1 = -1.
\end{align*}
The equation holds because the planned trajectory must avoid the huge cost at time step $2$. Solving this gives $u_0 = -\frac{1}{3}$. Thus, implementing MPC incurs a total cost that is at least $2u_0^2 = \frac{2}{9}$. In contrast, if one just pick $u_0 = 0$, the agent can pick $u_1$ based on the prediction revealed at time step $2$:
\begin{align*}
    u_1 = \begin{cases}
        0 & \text{ if }W_1 = 0,\\
        -1 & \text{ otherwise.}
    \end{cases}
\end{align*}
In this case, the expected cost incurred is $p$. Thus, we can claim that MPC is not the optimal policy when $p < \frac{2}{9}$. The underlying reason that MPC is suboptimal is because it does not consider what information may be available when we make the decision in the future. In this specific example, since $W_1$ is revealed at time $1$, we don't need to verify about the small probability event that leads to a huge loss.

We dive deeper into the reason why MPC \eqref{equ:MPC-in-expectation} is optimal in the LQR setting (Section \ref{sec:LQR}). Note that the expected optimal cost-to-go function at time step $1$ is
\begin{align}\label{MPC-suboptimal-example:e1}
    \Condexp{C_1^{\pi^\theta}(x; \Xi)}{I_1(\theta)} = \min_{u_1} \mathbb{E}\left[h_1(x, u_1) + h_2(X_2)\mid I_1(\theta)\right], \text{ s.t. } X_2 = x + u_1 + W_1.
\end{align}
Here, $u_1$ is $\mathcal{F}_1(\theta)$-measurable. And the true optimal policy at time $0$ is decided by solving
\begin{align*}
    \min_{u_0} h_0(x, u_0) + \mathbb{E}\left[C_1^{\pi^\theta}(X_1; \Xi)\mid I_0(\theta)\right], \text{ s.t. } X_1 = x + u_0.
\end{align*}
In general, we cannot use
\begin{align}\label{MPC-suboptimal-example:e2}
     \min_{u_1} \mathbb{E}\left[h_1(X_1, u_1) + h_2(X_2)\mid I_0(\theta)\right], \text{ s.t. } X_2 = X_1 + u_1 + W_1,
\end{align}
to replace $\Condexp{C_1^{\pi^\theta}(X_1; \Xi)}{I_0(\theta)}$ like what MPC does in \eqref{MPC-suboptimal-example:e0} because here $u_1$ is $\mathcal{F}_0(\theta)$-measurable in \eqref{MPC-suboptimal-example:e2}. Recall that $u_1$ is $\mathcal{F}_1(\theta)$-measurable in \eqref{MPC-suboptimal-example:e1} and $\mathcal{F}_0(\theta)$ is a subset of $\mathcal{F}_1(\theta)$. However, in the LQR setting, as the closed-form expression \eqref{thm:LQR_closed_form:e1}, the part of $\mathbb{E}\left[C_1^{\pi^\theta}(X_1; \Xi)\mid I_0(\theta)\right]$ that depends on $X_1$ will not change even if $\mathcal{F}_1(\theta)$ changes. Thus, we can assume $\mathcal{F}_1(\theta) = \mathcal{F}_0(\theta)$ without affecting the optimal action at time $0$. Therefore, MPC's replacement of $\mathbb{E}\left[C_1^{\pi^\theta}(X_1; \Xi)\mid I_0(\theta)\right]$ with \eqref{MPC-suboptimal-example:e2} is valid in the LQR setting.
\subsection{Infimal Convolution Properties}\label{appendix:infimal-conv-properties}
The first result states that the variant of infimal convolution preserves the strong convexity/smoothness of the input functions. The proof can be found in Appendix \ref{appendix:lemma:infimal-conv-variant}.

\begin{lemma}\label{lemma:infimal-conv-variant}
Consider a variant of infimal convolution defined as
\begin{align}\label{equ:infimal-conv-variant}
    \left(f \square_B \omega\right)(x) = \min_u \left\{f(u) + \omega(x - B u)\right\},
\end{align}
where $f: \mathbb{R}^m \to \mathbb{R}$, $\omega: \mathbb{R}^n \to \mathbb{R}$, and $B \in \mathbb{R}^{n \times m}$ is a matrix. Suppose that $f$ is a $\mu_f$-strongly convex function, and $\omega$ is a $\mu_\omega$-strongly convex and $\ell_\omega$-smooth function. Then, $f \square_B \omega$ is a $\left(\frac{\mu_\omega \mu_f}{\mu_f + \norm{B}^2 \mu_\omega}\right)$-strongly convex and $\ell_\omega$-smooth function. We also have $\nabla (f \square_B \omega)(x) = \nabla \omega(x - B u(x)).$
\end{lemma}

The second result is about the optimal solution of the variant of infimal convolution. It states that for some distributions, the covariance on the input will induce a variance on the optimal solution. We state it in Lemma \ref{lemma:infimal-conv-opt-solution-perturb-lower} and defer the proof to Appendix \ref{appendix:lemma:infimal-conv-pass-variance}.

\begin{lemma}\label{lemma:infimal-conv-opt-solution-perturb-lower}
Let $u_{(f\square_B \omega)}(x)$ denote the solution to the optimization problem \eqref{equ:infimal-conv-variant-def}. Suppose function $f$ is $\mu_f$-strongly convex. Function $\omega$ is $\mu_\omega$-strongly convex and $\ell_\omega$-smooth. Suppose $X$ is a random vector with bounded mean and $\Cov{X} = \Sigma \succeq \sigma_0 I$. Further, there exists a constant $C > 0$ such that for any positive integer $N$, $X$ can be decomposed as $X = \sum_{i=1}^N X_i$ for i.i.d. random vectors $X_i$ that satisfies $\mathbb{E}\left[\norm{X_i}^4\right] \leq C \cdot N^{-2}$. Then,
\begin{align*}
    \Tr{\Cov{u_{(f\square_B \omega)}(X)}} \geq \frac{n \sigma_0 \mu_\omega^2 \cdot \sigma_\text{min}(B)^2}{2(\ell_f + \ell_\omega \norm{B})^2}.
\end{align*}
\end{lemma}

As a remark, examples of $X$ that satisfies the assumptions include:
\begin{itemize}[leftmargin=*]
    \item Normal distribution $X \sim N(0, \Sigma)$. We have $X_i \sim N(0, \Sigma/N)$, thus $\mathbb{E}\left[\norm{X_i}^4\right] \leq 3\Tr{\Sigma} N^{-2}$.
    \item Poisson distribution (1D) with parameter $a$. We have $\Var{X} = a$ and $X_i$ follows Poisson distribution with parameter $a/N$. Thus, $\bE{X_i^4} = a^4 N^{-4}$.
\end{itemize}

The next result (Lemma \ref{lemma:infimal-conv-variant-grad-var-to-solution-var}) considers the case when there is an additional input $w$ to function $\omega$ in the infimal convolution. When this additional parameter causes a covariance on the gradient $\nabla_1 \omega(x, W)$, the optimal solution of the infimal convolution will also have a nonzero variance.

\begin{lemma}\label{lemma:infimal-conv-variant-grad-var-to-solution-var}
Suppose that $\omega(x, w)$ satisfies that $\omega(\cdot, w)$ is an $\ell_\omega$-smooth convex function for all $w$. For a random variable $W$, suppose that the following inequality holds for arbitrary fixed vector $x \in \mathbb{R}^n$,
\begin{align*}
    \Cov{\nabla_1\omega(x, W)} \succeq \sigma_0 I.
\end{align*}
Suppose that $f: \mathbb{R}^m \to \mathbb{R}$ is a $\mu_f$-strongly convex and $\ell_f$-smooth function ($m \leq n$). Let $B$ be a matrix in $\mathbb{R}^{n \times m}$. Then, the optimal solution of the infimal convolution 
\begin{align*}
    u_{(f\square_B \omega)}(x, w) \coloneqq \argmin_u\left(f(u) + \omega(x - B u, w)\right)
\end{align*}
satisfies that
\begin{align*}
    \Tr{\Cov{u_{(f\square_B \omega)}(x, W)}} \geq \frac{n \sigma_0 \cdot \sigma_\text{min}(B)^2}{2(\ell_f + \ell_\omega \norm{B})^2}.
\end{align*}
holds for arbitrary fixed vector $x$, where $\sigma_\text{min}(B)$ denotes the minimum singular value of $B$.
\end{lemma}

Lemma \ref{lemma:infimal-conv-variant-grad-var-to-solution-var} is useful for showing Lemma \ref{lemma:infimal-conv-opt-solution-perturb-lower}. We defer its proof to Appendix \ref{appendix:lemma:infimal-conv-variant-grad-var-to-solution-var}.
\subsection{Proof of Lemma \ref{lemma:well-condition-assump-verification-strongly-convex}}\label{appendix:lemma:well-condition-assump-verification-strongly-convex}
We use induction to show that $\mathbb{E}\left[C_t^{\pi^\theta}(x; \Xi)\mid I_t(\theta) = \iota_t(\theta)\right]$ is a $\mu_t$-strongly convex and $\ell_t$-smooth function for any $\iota_t(\theta)$, where the coefficients $\mu_t$ and $\ell_t$ are defined recursively in \eqref{lemma:well-condition-assump-verification-strongly-convex:recursive-def}. To simplify the notation, we will omit ``$I_t(\theta) =$'' in the conditional expectations throughout this proof when conditioning on a realization of the history $\iota_t(\theta)$.

Note that the statement holds for $t = T$, because $\mathbb{E}\left[C_T^{\pi^\theta}(x; \Xi)\mid \iota_T(\theta)\right] = h_T^x(x)$ and the terminal cost $h_T^x$ is $\mu_x$-strongly convex and $\ell_x$-smooth.

Suppose the statement holds for $t+1$. We see that
\begin{align*}
    \mathbb{E}\left[C_t^{\pi^\theta}(x; \Xi)\mid \iota_t(\theta)\right] = h_t^x(x) + \min_u\left(h_t^u(u) + \mathbb{E}\left[C_{t+1}^{\pi^\theta}(A_t x + B_t u + W_t; \Xi)\mid \iota_t(\theta)\right]\right).
\end{align*}
By the induction assumption, we know that $\mathbb{E}\left[C_{t+1}^{\pi^\theta}(\cdot; \Xi)\mid \iota_{t+1}(\theta)\right]$ is a $\mu_{t+1}$-strongly convex and $\ell_{t+1}$-smooth function for any $\iota_{t+1}(\theta)$. Thus, $\mathbb{E}\left[C_{t+1}^{\pi^\theta}(\cdot + W_t; \Xi)\mid \iota_{t}(\theta)\right]$ is also a $\mu_{t+1}$-strongly convex and $\ell_{t+1}$-smooth function. Therefore,
\begin{align*}
    \min_u\left(h_t^u(u) + \mathbb{E}\left[C_{t+1}^{\pi^\theta}(x + B_t u + W_t; \Xi)\mid \iota_t(\theta)\right]\right)
\end{align*}
is a $\frac{\mu_u \mu_{t+1}}{\mu_u + b^2 \mu_{t+1}}$-strongly convex and $\ell_{t+1}$-smooth function of $x$ by Lemma \ref{lemma:infimal-conv-variant}. By changing the variable from $x$ to $A_t x$, we see that
\begin{align*}
    \min_u\left(h_t^u(u) + \mathbb{E}\left[C_{t+1}^{\pi^\theta}(A_t x + B_t u + W_t; \Xi)\mid \iota_t(\theta)\right]\right)
\end{align*}
is a $\mu_A \cdot \frac{\mu_u \mu_{t+1}}{\mu_u + b^2 \mu_{t+1}}$-strongly convex and $\ell_A \cdot \ell_{t+1}$-smooth function by Assumption \ref{assump:well-conditioned-cost}. Since $h_t^x$ is a $\mu_x$-strongly convex and $\ell_x$-smooth function, we see that $\mathbb{E}\left[C_t^{\pi^\theta}(x; \Xi)\mid \iota_t(\theta)\right]$ is also a $\mu_t$-strongly convex and $\ell_t$-smooth function because
\begin{align*}
    \mu_t = \mu_x + \mu_A \cdot \frac{\mu_u \mu_{t+1}}{\mu_u + b^2 \mu_{t+1}}, \text{ and } \ell_t = \ell_x + \ell_A \cdot \ell_{t+1}.
\end{align*}
\subsection{Proof of Theorem \ref{coro:one-step-pred-power-lower}}\label{appendix:coro:one-step-pred-power-lower}
Note that the optimal action at time step $t$ is determined by
\begin{align}\label{equ:joint-gaussian-intuition:e1}
    \pi_t^\theta(x; I_t(\theta)) \coloneqq \argmin_u \left(h_t^u(u) + \mathbb{E}\left[C_{t+1}^{\pi^\theta}(A_t x + B_t u + W_t; \Xi)\mid I_t(\theta)\right]\right).
\end{align}
This can be further simplified to
\begin{align*}
    \pi_t^\theta(x; I_t(\theta)) &\coloneqq \argmin_u \left(h_t^u(u) + \bar{C}_{t+1}^{\pi^\theta}(A_t x + B_t u + W_{t\mid t}^\theta)\right).
\end{align*}
The additional input $I_t(\theta)$ is not required for $\bar{C}_{t+1}^{\pi^\theta}$ because the function $\bar{C}_{t+1}^{\pi^\theta}(x; \iota_t(\theta))$ does not change with the history $\iota_t(\theta)$ under Assumption \ref{assump:one-step-prediction-indep-pairs}. The reason is that $W_t - W_{t\mid t}^\theta$ and all future predictions and disturbances $W_{t+1:T-1}, V_{t+1:T-1}^\theta$ are independent with the history $I_t(\theta)$.

By \eqref{equ:joint-gaussian-intuition:e2}, we see that
\begin{align*}
    \pi_t^\theta(x; I_t(\theta)) = u_{(h_t^u\square_{-B_t}\bar{C}_{t+1}^{\pi^\theta})}(A_t x + W_{t\mid t}^\theta).
\end{align*}
Under Assumption \ref{assump:one-step-prediction-indep-pairs}, we see that
\begin{align*}
    \Cov{W_{t\mid t}^\theta} = \Cov{W_t} - \Cov{W_t\mid V_t(\theta)}  \succeq \lambda_t(\theta) I
\end{align*}
and $W_{t\mid t}^\theta$ is Gaussian. Therefore, we can apply Lemma \ref{lemma:infimal-conv-opt-solution-perturb-lower} to obtain that
\begin{align*}
    \Tr{\Cov{\pi_t^\theta(x; I_t(\theta))\mid \mathcal{F}_t(0)}} \geq \sigma_t \coloneqq \frac{n \lambda_t(\theta) \mu_{t+1}^2 \cdot \mu_B}{2(\ell_u + \ell_{t+1}\sqrt{\ell_B})^2}.
\end{align*}
Thus, Condition \ref{cond:optimal-policy-psd-covariance}\,\ref{cond:optimal-policy-psd-covariance:b} holds with $\sigma_t$.

On the other hand, Condition \ref{cond:Q-function-quadratic-growth} holds with $M_t = \mu_t I$ by Lemma \ref{lemma:well-condition-assump-verification-strongly-convex}. Therefore, by Theorem \ref{thm:Bellman-improvement-not-worse}, we obtain that $P(\theta) \geq \sum_{t=0}^{T-1} \mu_u \sigma_t$.
\subsection{Proof of Lemma \ref{lemma:infimal-conv-variant}}\label{appendix:lemma:infimal-conv-variant}
By the definition of conjugate, we see that
\begin{subequations}\label{lemma:infimal-conv-variant:e1}
\begin{align}
    \left(f \square_B \omega\right)^*(y)
    ={}& \max_x\left\{\langle y, x\rangle - \min_u\left\{f(u) + \omega(x - B u)\right\}\right\}\label{lemma:infimal-conv-variant:e1:s1}\\
    ={}& \max_x \max_u\left\{\langle y, x\rangle - f(u) - \omega(x - B u)\right\}\nonumber\\
    ={}& \max_x \max_u\left\{\langle y, x - Bu\rangle + \langle y, B u\rangle - f(u) - \omega(x - B u)\right\}\nonumber\\
    ={}& \max_u \max_x\left\{\left(\langle y, x - Bu\rangle - \omega(x - B u)\right) + \left(\langle B^\top y, u\rangle - f(u)\right)\right\}\label{lemma:infimal-conv-variant:e1:s2}\\
    ={}& \max_u \left\{\max_x\left\{\langle y, x - Bu\rangle - \omega(x - B u)\right\} + \langle B^\top y, u\rangle - f(u)\right\}\nonumber\\
    ={}& \max_u\left\{\omega^*(y) + \langle B^\top y, u\rangle - f(u)\right\}\label{lemma:infimal-conv-variant:e1:s3}\\
    ={}& \omega^*(y) + f^*(B^\top y),\label{lemma:infimal-conv-variant:e1:s4}
\end{align} 
\end{subequations}
where we use the definition of $f \square_B \omega$ in \eqref{lemma:infimal-conv-variant:e1:s1}; we change the order of taking the maximum and use $\langle y, B u\rangle = \langle B^\top y, u\rangle$ in \eqref{lemma:infimal-conv-variant:e1:s2}; we use the definition of $\omega^*$ in \eqref{lemma:infimal-conv-variant:e1:s3}; we use the definition of $f^*$ in \eqref{lemma:infimal-conv-variant:e1:s4}.

Since $f \square_B \omega$ is convex, by Theorem 4.8 in \cite{beck2017first}, we know that
\begin{align}\label{lemma:infimal-conv-variant:e2}
    \left(f \square_B \omega\right)(y) = \left(\omega^*(y) + f^*(B^\top y)\right)^*.
\end{align}
Since $\omega$ is a $\mu_\omega$-strongly convex and $\ell_\omega$-smooth function, we know $\omega^*$ is an $\frac{1}{\ell_\omega}$-strongly convex and $\frac{1}{\mu_\omega}$-smooth function by the conjugate correspondence theorem \citep{beck2017first}. Similarly, we know that $f^*$ is a $\frac{1}{\mu_f}$-smooth convex function. Thus, we know that $\omega^*(y) + f^*(B^\top y)$ is an $\frac{1}{\ell_\omega}$-strongly convex and $\left(\frac{1}{\mu_\omega} + \frac{\norm{B}^2}{\mu_f}\right)$-smooth function. Therefore, by the conjugate correspondence theorem, we know that $f \square_B \omega$ is a $\left(\frac{\mu_\omega \mu_f}{\mu_f + \norm{B}^2 \mu_\omega}\right)$-strongly convex and $\ell_\omega$-smooth function.

Now, we show that
\begin{align}\label{lemma:infimal-conv-variant:e3}
    \nabla (f \square_B \omega)(x) = \nabla \omega(x - B u(x)).
\end{align}
Following a similar approach with the proof of Theorem 5.30 in \cite{beck2017first}, we define $z = \nabla \omega(x - B u(x))$. Define function $\phi(\xi) \coloneqq (f \square_B \omega)(x+\xi) - (f \square_B \omega)(x) - \langle \xi, z\rangle$. We see that
\begin{subequations}\label{lemma:infimal-conv-variant:e4}
\begin{align}
    \phi(\xi) ={}& (f\square_B \omega)(x + \xi) - (f\square_B \omega)(x) - \langle \xi, z\rangle\nonumber\\
    \leq{}& \omega(x + \xi - B u(x)) - \omega(x - B u(x)) - \langle \xi, z\rangle \label{lemma:infimal-conv-variant:e4:s1}\\
    \leq{}& \langle \xi, \nabla \omega(x + \xi - B u(x))\rangle - \langle \xi, z\rangle \label{lemma:infimal-conv-variant:e4:s2}\\
    ={}& \langle \xi, \nabla \omega(x + \xi - B u(x)) - \nabla \omega(x - B u(x))\rangle \nonumber\\
    \leq{}& \norm{\xi} \cdot \norm{\nabla \omega(x + \xi - B u(x)) - \nabla \omega(x - B u(x))} \label{lemma:infimal-conv-variant:e4:s3}\\
    \leq{}& \ell_\omega \norm{\xi}^2, \label{lemma:infimal-conv-variant:e4:s4}
\end{align}
\end{subequations}
where in \eqref{lemma:infimal-conv-variant:e4:s1}, we use
\begin{align*}
    (f\square_B \omega)(x + \xi) \leq{}& f(u(x)) + \omega(x + \xi - B u(x)), \text{ and }\\
    (f\square_B \omega)(x) ={}& f(u(x)) + \omega(x - B u(x));
\end{align*}
we use the convexity of $\omega$ in \eqref{lemma:infimal-conv-variant:e4:s2}; we use the Cauchy-Schwarz inequality in \eqref{lemma:infimal-conv-variant:e4:s3}; we use the assumption that $\omega$ is $\ell_\omega$-smooth in \eqref{lemma:infimal-conv-variant:e4:s4}.

Since $(f \square_B \omega)$ is a convex function, $\phi$ is also convex, thus we see that
\begin{align*}
    \phi(\xi) \geq 2\phi(0) - \phi(-\xi) = - \phi(-\xi) \geq - \ell_\omega \norm{\xi}^2.
\end{align*}
Combining this with \eqref{lemma:infimal-conv-variant:e4}, we conclude that $\lim_{\norm{\xi} \to 0}\abs{\phi(\xi)}/\norm{\xi} = 0$. Thus, \eqref{lemma:infimal-conv-variant:e3} holds.
\subsection{Proof of Lemma \ref{lemma:infimal-conv-opt-solution-perturb-lower}}\label{appendix:lemma:infimal-conv-pass-variance}
By Theorem \ref{thm:variance-pass-through-function}, we see that
\begin{align*}
    \Cov{\nabla \omega(X)} \geq \sigma_0 \mu_\omega^2.
\end{align*}
Then, we apply Lemma \ref{lemma:infimal-conv-variant-grad-var-to-solution-var} with the second function input to the infimal convolution as $\tilde{\omega}(x, w) \coloneqq \omega(x + w)$. In the context of Lemma \ref{lemma:infimal-conv-variant-grad-var-to-solution-var}, we set $W = X$, so the assumption about the covariance of the gradient holds with
\begin{align*}
    \Cov{\nabla_1 \tilde{\omega}(x, W)} \succeq \sigma_0 \mu_\omega^2.
\end{align*}
Note that for any fixed $w$, $\tilde{\omega}(\cdot, w)$ is $\mu_\omega$-strongly convex. Therefore, we obtain that
\begin{align*}
    \Tr{\Cov{u_{(f\square_B \omega)}(X)}} = \Tr{\Cov{u_{(f\square_B \tilde{\omega})}(0, W)}} \geq \frac{n \sigma_0 \mu_\omega^2 \cdot \sigma_\text{min}(B)^2}{2(\ell_f + \ell_\omega \norm{B})^2}
\end{align*}
\subsection{Proof of Lemma \ref{lemma:infimal-conv-variant-grad-var-to-solution-var}}\label{appendix:lemma:infimal-conv-variant-grad-var-to-solution-var}
Because function $c$ is $\ell_c$-smooth, we have
\begin{align}\label{lemma:infimal-conv-variant-grad-var-to-solution-var:e1}
    \norm{\nabla c(u(x, w)) - \nabla c(u(x, w'))} \leq \ell_c \norm{u(x, w) - u(x, w')}
\end{align}
Because function $f$ is $\ell_f$-smooth, we have
\begin{subequations}\label{lemma:infimal-conv-variant-grad-var-to-solution-var:e2}
\begin{align}
    &\norm{B^\top \nabla_1 f(x - B u(x, w), w) - B^\top \nabla_1 f(x - B u(x, w'), w')}\nonumber\\
    \geq{}& \norm{B^\top \nabla_1 f(x - B \cdot \mathbb{E}_W\left[u(x, W)\right], w) - B^\top \nabla_1 f(x - B \cdot \mathbb{E}_W\left[u(x, W)\right], w')}\nonumber\\
    &- \norm{B^\top \nabla_1 f(x - B \cdot u(x, w), w) - B^\top \nabla_1 f(x - B \cdot \mathbb{E}_W\left[u(x, W)\right], w)}\nonumber\\
    &- \norm{B^\top \nabla_1 f(x - B \cdot u(x, w'), w') - B^\top \nabla_1 f(x - B \cdot \mathbb{E}_W\left[u(x, W)\right], w')}\label{lemma:infimal-conv-variant-grad-var-to-solution-var:e2:s1}\\
    \geq{}& \norm{B^\top \nabla_1 f(x - B \cdot \mathbb{E}_W\left[u(x, W)\right], w) - B^\top \nabla_1 f(x - B \cdot \mathbb{E}_W\left[u(x, W)\right], w')}\nonumber\\
    & - \ell_f \norm{B} \cdot \left(\norm{u(x, w) - \mathbb{E}_W\left[u(x, W)\right]} + \norm{u(x, w') - \mathbb{E}_W\left[u(x, W)\right]}\right),\label{lemma:infimal-conv-variant-grad-var-to-solution-var:e2:s2}
\end{align}
\end{subequations}
where we use the triangle inequality in \eqref{lemma:infimal-conv-variant-grad-var-to-solution-var:e2:s1}; we use the smoothness of $f$ in \eqref{lemma:infimal-conv-variant-grad-var-to-solution-var:e2:s2}.

Note that by the first-order optimality condition, we have
\begin{align*}
    \nabla c(u(x, w)) - B^\top \nabla_1 f(x - B\cdot u(x, w), w) = 0.
\end{align*}
Therefore, for any $w, w'$, we have that
\begin{align}\label{lemma:infimal-conv-variant-grad-var-to-solution-var:e3}
    \nabla c(u(x, w)) - \nabla c(u(x, w')) = B^\top \nabla_1 f(x - B \cdot u(x, w), w) - B^\top \nabla_1 f(x - B \cdot u(x, w'), w').
\end{align}
By combining \eqref{lemma:infimal-conv-variant-grad-var-to-solution-var:e3} with \eqref{lemma:infimal-conv-variant-grad-var-to-solution-var:e1} and \eqref{lemma:infimal-conv-variant-grad-var-to-solution-var:e2}, we obtain that
\begin{align*}
    &\ell_c \norm{u(x, w) - u(x, w')}\\
    &+ \ell_f \cdot \norm{B} \cdot \left(\norm{u(x, w) - \mathbb{E}_W\left[u(x, W)\right]} + \norm{u(x, w') - \mathbb{E}_W\left[u(x, W)\right]}\right)\\
    \geq{}& \norm{B^\top \nabla_1 f(x - B \cdot \mathbb{E}_W\left[u(x, W)\right], w) - B^\top \nabla_1 f(x - B \cdot \mathbb{E}_W\left[u(x, W)\right], w')}
\end{align*}
holds for arbitrary $w$ and $w'$.
Let $W'$ be a random vector independent of $W$ and have the same distribution. By replacing $w/w'$ with $W/W'$ respectively, we see
\begin{align*}
    &\ell_c \norm{u(x, W) - u(x, W')}\\
    &+ \ell_f \cdot \norm{B} \cdot \left(\norm{u(x, W) - \mathbb{E}_W\left[u(x, W)\right]} + \norm{u(x, W') - \mathbb{E}_W\left[u(x, W)\right]}\right)\\
    \geq{}& \norm{B^\top \nabla_1 f(x - B \cdot \mathbb{E}_W\left[u(x, W)\right], W) - B^\top \nabla_1 f(x - B \cdot \mathbb{E}_W\left[u(x, W)\right], W')},
\end{align*}
which implies
\begin{align}\label{lemma:infimal-conv-variant-grad-var-to-solution-var:e4}
    &(\ell_c + \ell_f\norm{B}) \left(\norm{u(x, W) - \mathbb{E}_W\left[u(x, W)\right]} + \norm{u(x, W') - \mathbb{E}_W\left[u(x, W)\right]}\right)\nonumber\\
    \geq{}& \norm{B^\top \nabla_1 f(x - B \cdot \mathbb{E}_W\left[u(x, W)\right], W) - B^\top\nabla_1 f(x - B \cdot \mathbb{E}_W\left[u(x, W)\right], W')}
\end{align}
by the triangle inequality. Taking the square of both sides of \eqref{lemma:infimal-conv-variant-grad-var-to-solution-var:e4} and applying the AM-GM inequality gives that
\begin{align}\label{lemma:infimal-conv-variant-grad-var-to-solution-var:e5}
    &2(\ell_c + \ell_f\norm{B})^2 \norm{u(x, W) - \mathbb{E}_W\left[u(x, W)\right]}^2 + 2(\ell_c + \ell_f\norm{B})^2 \norm{u(x, W') - \mathbb{E}_W\left[u(x, W)\right]}^2\nonumber\\
    \geq{}& \norm{B^\top \nabla_1 f(x - B \cdot \mathbb{E}_W\left[u(x, W)\right], W) - B^\top \nabla_1 f(x - B \cdot \mathbb{E}_W\left[u(x, W)\right], W')}^2.
\end{align}
Let $Y \coloneqq \nabla_1 f(x - B \cdot \mathbb{E}_W\left[u(x, W)\right], W) - \nabla_1 f(x - B \cdot \mathbb{E}_W\left[u(x, W)\right], W')$. Note that the right-hand side of \eqref{lemma:infimal-conv-variant-grad-var-to-solution-var:e5} can be expressed as $\norm{B^\top Y}^2 = \Tr{B^\top (Y Y^\top) B}$. By taking the expectations of both sides, we obtain that
\begin{align*}
    4 (\ell_c + \ell_f\norm{B})^2 \Tr{\Cov{u(x, W)}} \geq{}& 2 \Tr{B^\top \Cov{\nabla_1 f(x - B \mathbb{E}_W\left[u(x, W)\right], W)} B}\\
    \geq{}& 2 n \sigma_0 \sigma_\text{min}(B)^2.
\end{align*}
In the last inequality, we use the property that the trace of a positive semi-definite matrix equals the sum of its eigenvalues. Thus, it is greater than or equal to $n$ times the smallest eigenvalue $\sigma_0 \sigma_\text{min}(B)^2$. Rearranging the terms finishes the proof.

\section{Useful Technical Results}\label{appendix:useful-lemmas}
In this section, we state a useful result about what functions can pass the covariance of its input to the output in Theorem \ref{thm:variance-pass-through-function}, which is used to show \Cref{lemma:infimal-conv-opt-solution-perturb-lower}. We defer the proof to Appendix \ref{appendix:thm:variance-pass-through-function}.

\begin{theorem}\label{thm:variance-pass-through-function}
Suppose that a function $g:\mathbb{R}^d\to \mathbb{R}^d$ satisfies
\begin{align}\label{equ:inequality-question-01-03}
    \langle g(x) - g(x'), x - x'\rangle \geq \gamma \norm{x - x'}^2, \text{ and } \norm{g(x) - g(x')} \leq L \norm{x - x'},\;\forall\,x,x'\in\mathbb{R}^d.
\end{align}
Additionally, there exists a positive constant $\ell$ such that
\begin{align}\label{equ:inequality-question-01-03:e0}
    - \ell I \preceq \nabla^2 g_i(x) \preceq \ell I, \;\forall\,x\in\mathbb{R}^d, i \in [d].
\end{align}
Suppose $X$ is a random vector that satisfies $\abs{\bE{X}} < \infty$ and $\Cov{X} = \Sigma \succeq \mu I$. Further, there exists a constant $C > 0$ such that for any positive integer $N$, $X$ can be decomposed as $X = \sum_{i=1}^N X_i$ for i.i.d. random vectors $X_i$ that satisfies $\mathbb{E}\left[\norm{X_i}^4\right] \leq C \cdot N^{-2}$. Then, we have 
\begin{align*}
    \Cov{g(X)} \succeq \mu \gamma^2 I.
\end{align*}
\end{theorem}
As a remark, the gradient of a well-conditioned function satisfies the conditions in \eqref{equ:inequality-question-01-03}.
\subsection{Proof of Theorem \ref{thm:variance-pass-through-function}}\label{appendix:thm:variance-pass-through-function}
Without any loss of generality, we assume $\bE{X} = 0$ because we can view $g(\bE{X} + \cdot)$ as the function and subtract the mean from the random variables. The assumptions about $g$ and $X$ in Theorem \ref{thm:variance-pass-through-function} still hold.

For any $i\in [d]$ and $\epsilon\in\mathbb{R}^d$, we have the Taylor series expansion Lagrangian form (see Chapter 3.2 of \cite{marsden2003vector}) %\zc{replace this citation by a math textbook}
\begin{align}\label{equ:inequality-question-01-03:e1}
    g_i(x + \epsilon) = g_i(x) + \nabla g_i(x)^\top \epsilon + \frac{1}{2} \epsilon^\top \nabla^2 g_i(\bar{x}^{(i)}) \epsilon,
\end{align}
where $\bar{x}^{(i)}$ is a point on the line segment between $x$ and $x+\epsilon$. For notational convenience, let
\begin{align*}
    \nabla g(x) \coloneqq \begin{bmatrix}
        \nabla g_1(x)^\top\\
        \vdots\\
        \nabla g_d(x)^\top
    \end{bmatrix} \in \mathbb{R}^{d\times d}, \text{ and } v_1(x, \epsilon) \coloneqq \begin{bmatrix}
        \epsilon^\top \nabla^2 g_1(\bar{x}^{(1)}) \epsilon\\
        \vdots\\
        \epsilon^\top \nabla^2 g_d(\bar{x}^{(d)}) \epsilon
    \end{bmatrix} \in \mathbb{R}^d.
\end{align*}
With the above notation, Eq. \eqref{equ:inequality-question-01-03:e1} can be equivalently written as
\begin{align}\label{equ:inequality-question-01-03:e2}
    g(x+\epsilon) - g(x) = \nabla g(x) \cdot \epsilon + \frac{1}{2}v_1(x, \epsilon).
\end{align}
From Eq. \eqref{equ:inequality-question-01-03:e0}, we know that $\abs{v(x, \epsilon)_i} \leq \ell \norm{\epsilon}^2$, which implies
\begin{align}\label{equ:inequality-question-01-03:e2-1}
    \norm{v_1(x, \epsilon)} \leq \ell\sqrt{d} \norm{\epsilon}^2.
\end{align}
In addition, by Eq. \eqref{equ:inequality-question-01-03}, we see that
\begin{align*}
    \langle g(x + \epsilon) - g(x), \epsilon\rangle \geq \gamma \norm{\epsilon}^2.
\end{align*}
Substituting Eq. \eqref{equ:inequality-question-01-03:e2} into the above equation and rearranging the terms, we obtain
\begin{align*}
    \epsilon^\top \cdot \nabla g(x) \cdot \epsilon \geq \gamma \norm{\epsilon}^2 - \epsilon^\top \cdot v_1(x, \epsilon),
\end{align*}
which is equivalent to
\begin{align*}
    \epsilon^\top \cdot \frac{\nabla g(x) + \nabla g(x)^\top}{2} \cdot \epsilon \geq \gamma \norm{\epsilon}^2 - \epsilon^\top \cdot v_1(x, \epsilon).
\end{align*}
Observe that the term subtracted from the right-hand side satisfies $\abs{\epsilon^\top \cdot v_1(x, \epsilon)} \leq \ell \sqrt{d} \norm{\epsilon}^3$, which follows from Cauchy–Schwarz inequality and Eq. (\ref{equ:inequality-question-01-03:e2-1}). Therefore, since the previous inequality holds for any $\epsilon\in\mathbb{R}^d$, taking $\epsilon \to 0$ gives that
\begin{align}\label{equ:inequality-question-01-03:e3}
    \frac{\nabla g(x) + \nabla g(x)^\top}{2} \succeq \gamma I.
\end{align}
Before we proceed, we first state and prove a lemma that can convert the summation in Eq. \eqref{equ:inequality-question-01-03:e3} into a product form. 
\begin{lemma}\label{lemma:matrix-pd-sum-to-product}
Let $M \in \mathbb{R}^{d \times d}$ be a real-valued matrix satisfying $M + M^\top \succeq 2 \gamma I$. Then, for any positive definite matrix $\Sigma \succeq \mu I$, we have $M \Sigma M^\top \succeq \mu \gamma^2 I$.
\end{lemma}
\begin{proof}[Proof of Lemma \ref{lemma:matrix-pd-sum-to-product}]
Since $M + M^\top \succeq 2 \gamma I$, we have for any $x\in\mathbb{R}^d$ that
\begin{align*}
    2\gamma\|x\|^2&\leq 2x^\top M^\top x=2x^\top \Sigma^{-1/2}\Sigma^{1/2}M^\top x\leq 2\|\Sigma^{-1/2}x\|\|\Sigma^{1/2}M^\top x\|\\
    &\leq 2\mu^{-1/2}\|x\|\|\Sigma^{1/2}M^\top x\|,
\end{align*}
where the last inequality follows from
$\Sigma \succeq \mu I \;\Rightarrow\;\|\Sigma^{-1/2}x\|=\sqrt{x^\top \Sigma^{-1} x}\leq \mu^{-1/2}\|x\|$. Rearranging terms, we obtain
\begin{align*}
    \gamma\mu^{1/2}\|x\|\leq \|\Sigma^{1/2}M^\top x\|.
\end{align*}
Squaring both sides concludes the proof.
\end{proof}

Next, we state and prove a lemma about the lower bound of the covariance induced by an additive random noise on the input that is useful when the noise is sufficiently small.

\begin{lemma}\label{lemma:small-noise-covariance-lower-bound}
Let $\varepsilon$ be a mean-zero random vector in $\mathbb{R}^d$ that satisfies $\underline{\delta} I \preceq \Cov{\varepsilon}$ and $\mathbb{E}\left[\norm{\varepsilon}^4\right] \leq \overline{\gamma}$. Let $g$ be a function that satisfies \eqref{equ:inequality-question-01-03} and \eqref{equ:inequality-question-01-03:e0}. Then, for arbitrary fixed real vector $x \in \mathbb{R}^d$, we have
\begin{align*}
    \Cov{g(x + \varepsilon)} \succeq \left(\gamma^2 \underline{\delta} - 2 L \ell d^2 \cdot \overline{\gamma}^{\frac{3}{4}} - \ell^2 d \overline{\gamma}\right) I.
\end{align*}
\end{lemma}
\begin{proof}[Proof of Lemma \ref{lemma:small-noise-covariance-lower-bound}]
We first derive bounds on the $i$ th moment of $\norm{\varepsilon}$ ($i = 1, 2, 3$). By Jensen's inequality, we have
\begin{align}\label{lemma:small-noise-covariance-lower-bound:e0-2}
    \mathbb{E}\left[\norm{\varepsilon}^2\right] = \mathbb{E}\left[\left(\norm{\varepsilon}^4\right)^{\frac{1}{2}}\right] \leq \left(\bE{\norm{\varepsilon}^4}\right)^{\frac{1}{2}} \leq \overline{\gamma}^{\frac{1}{2}}.
\end{align}
Using Jensen'e inequality again, we obtain that
\begin{align}\label{lemma:small-noise-covariance-lower-bound:e0-1}
    \mathbb{E}\left[\norm{\varepsilon}\right] \leq \left(\bE{\norm{\varepsilon}^2}\right)^{\frac{1}{2}} \leq \overline{\gamma}^{\frac{1}{4}}.
\end{align}
Lastly, by the Cauchy-Schwartz inequality, we see that
\begin{align}\label{lemma:small-noise-covariance-lower-bound:e0-3}
    \mathbb{E}\left[\norm{\varepsilon}^3\right] \leq \left(\mathbb{E}\left[\norm{\varepsilon}^4\right] \cdot \mathbb{E}\left[\norm{\varepsilon}^2\right]\right)^{\frac{1}{2}} \leq \overline{\gamma}^{\frac{3}{4}}.
\end{align}

Note that by \eqref{equ:inequality-question-01-03:e2}, we have
\begin{align}\label{lemma:small-noise-covariance-lower-bound:e1}
    \Cov{g(x + \varepsilon)} = \Cov{g(x + \varepsilon) - g(x)} = \Cov{\nabla g(x) \cdot \varepsilon + \frac{1}{2}v_1(x, \varepsilon)}.
\end{align}
Since $\mathbb{E}\left[\varepsilon\right] = 0$, we can further decompose \eqref{lemma:small-noise-covariance-lower-bound:e1} as
\begin{align}\label{lemma:small-noise-covariance-lower-bound:e2}
    &\Cov{\nabla g(x) \cdot \varepsilon + \frac{1}{2}v_1(x, \varepsilon)}\nonumber\\
    ={}& \mathbb{E}\left[\left(\nabla g(x) \cdot \varepsilon + \frac{1}{2}v_1(x, \varepsilon) - \frac{1}{2}\mathbb{E}\left[v_1(x, \varepsilon)\right]\right)\left(\nabla g(x) \cdot \varepsilon + \frac{1}{2}v_1(x, \varepsilon) - \frac{1}{2}\mathbb{E}\left[v_1(x, \varepsilon)\right]\right)^\top\right]\nonumber\\
    ={}& \nabla g(x) \cdot \Cov{\varepsilon} \cdot \nabla g(x)^\top + \nabla g(x)\cdot \mathbb{E}\left[\varepsilon \cdot \left(\frac{1}{2}v_1(x, \varepsilon) - \frac{1}{2}\mathbb{E}\left[v_1(x, \varepsilon)\right]\right)^\top\right]\nonumber\\
    &+ \mathbb{E}\left[\left(\frac{1}{2}v_1(x, \varepsilon) - \frac{1}{2}\mathbb{E}\left[v_1(x, \varepsilon)\right]\right)\cdot \varepsilon^\top\right] \cdot \nabla g(x)^\top + \frac{1}{4}\Cov{v_1(x, \varepsilon)}.
\end{align}
By Lemma \ref{lemma:matrix-pd-sum-to-product} and \eqref{equ:inequality-question-01-03:e3}, we know the first term in \eqref{lemma:small-noise-covariance-lower-bound:e2} is lower bounded by
\begin{align}
    \nabla g(x) \cdot \Cov{\varepsilon} \cdot \nabla g(x)^\top \succeq \gamma^2 \underline{\delta} I.
\end{align}
Define the residual term as the sum of the last 3 terms in \eqref{lemma:small-noise-covariance-lower-bound:e2}:
\begin{align}\label{lemma:small-noise-covariance-lower-bound:e3-0}
    R\coloneqq{}& \nabla g(x)\cdot \mathbb{E}\left[\varepsilon \cdot \left(\frac{1}{2}v_1(x, \varepsilon) - \frac{1}{2}\mathbb{E}\left[v_1(x, \varepsilon)\right]\right)^\top\right]\nonumber\\
    &+ \mathbb{E}\left[\left(\frac{1}{2}v_1(x, \varepsilon) - \frac{1}{2}\mathbb{E}\left[v_1(x, \varepsilon)\right]\right)\cdot \varepsilon^\top\right] \cdot \nabla g(x)^\top + \frac{1}{4}\Cov{v_1(x, \varepsilon)}.
\end{align}
To show Lemma \ref{lemma:small-noise-covariance-lower-bound}, we only need to show
\begin{align}\label{lemma:small-noise-covariance-lower-bound:e3-1}
    \norm{R} \leq 2 L \ell d^2 \cdot \overline{\gamma}^{\frac{3}{4}} + \ell^2 d \overline{\gamma}.
\end{align}
To see this, note that
\begin{subequations}\label{lemma:small-noise-covariance-lower-bound:e3}
\begin{align}
    &\norm{\nabla g(x)\cdot \mathbb{E}\left[\varepsilon \cdot \left(\frac{1}{2}v_1(x, \varepsilon) - \frac{1}{2}\mathbb{E}\left[v_1(x, \varepsilon)\right]\right)^\top\right]}\nonumber\\
    \leq{}& \norm{\nabla g(x)} \cdot \norm{\mathbb{E}\left[\varepsilon \cdot \left(\frac{1}{2}v_1(x, \varepsilon) - \frac{1}{2}\mathbb{E}\left[v_1(x, \varepsilon)\right]\right)^\top\right]}\label{lemma:small-noise-covariance-lower-bound:e3:s1}\\
    \leq{}& \frac{L}{2} \left(\norm{\mathbb{E}\left[\varepsilon \cdot v_1(x, \varepsilon)^\top\right]} + \norm{\mathbb{E}\left[\varepsilon\right] \cdot \mathbb{E}\left[v_1(x, \varepsilon)\right]^\top}\right)\label{lemma:small-noise-covariance-lower-bound:e3:s2}\\
    \leq{}& \frac{L}{2} \left(\mathbb{E}\left[\norm{\varepsilon \cdot v_1(x, \varepsilon)^\top}\right] + \norm{\mathbb{E}\left[\varepsilon\right]} \cdot \norm{\mathbb{E}\left[v_1(x, \varepsilon)\right]}\right)\label{lemma:small-noise-covariance-lower-bound:e3:s3}\\
    \leq{}& \frac{L}{2} \left(\mathbb{E}\left[\norm{\varepsilon} \cdot \norm{v_1(x, \varepsilon)}\right] + \mathbb{E}\left[\norm{\varepsilon}\right] \cdot \mathbb{E}\left[\norm{v_1(x, \varepsilon)}\right]\right)\label{lemma:small-noise-covariance-lower-bound:e3:s4}\\
    \leq{}& \frac{L \ell \sqrt{d}}{2} \cdot \left(\mathbb{E}\left[\norm{\varepsilon}^3\right] + \mathbb{E}\left[\norm{\varepsilon}\right] \cdot \mathbb{E}\left[\norm{\varepsilon}^2\right]\right)\label{lemma:small-noise-covariance-lower-bound:e3:s5}\\
    \leq{}& L \ell d^2 \cdot \overline{\gamma}^{\frac{3}{4}}, \label{lemma:small-noise-covariance-lower-bound:e3:s6}
\end{align}
\end{subequations}
where we use the definition of the induced matrix norm in \eqref{lemma:small-noise-covariance-lower-bound:e3:s1}; we use \eqref{equ:inequality-question-01-03} and the triangle inequality in \eqref{lemma:small-noise-covariance-lower-bound:e3:s2}; we use the Jensen's inequality and the definition of the induced matrix norm in \eqref{lemma:small-noise-covariance-lower-bound:e3:s3} and \eqref{lemma:small-noise-covariance-lower-bound:e3:s4}; we use \eqref{equ:inequality-question-01-03:e2-1} in \eqref{lemma:small-noise-covariance-lower-bound:e3:s5}; we use the bounds on the moments of $\norm{\varepsilon}$ \eqref{lemma:small-noise-covariance-lower-bound:e0-2}, \eqref{lemma:small-noise-covariance-lower-bound:e0-1}, and \eqref{lemma:small-noise-covariance-lower-bound:e0-3} in \eqref{lemma:small-noise-covariance-lower-bound:e3:s6}.

On the other hand, we know that $\Cov{v_1(x, \varepsilon)}$ is a positive semi-definite matrix that satisfies
\begin{align*}
    \Cov{v_1(x, \varepsilon)} = \mathbb{E}\left[v_1(x, \varepsilon) v_1(x, \varepsilon)^\top\right] - \mathbb{E}\left[v_1(x, \varepsilon)\right] \cdot \mathbb{E}\left[v_1(x, \varepsilon)\right]^\top \preceq \mathbb{E}\left[v_1(x, \varepsilon) v_1(x, \varepsilon)^\top\right].
\end{align*}
Thus, its induced matrix norm can be upper bounded by
\begin{align*}
    \norm{\Cov{v_1(x, \varepsilon)}} \leq{}& \norm{\mathbb{E}\left[v_1(x, \varepsilon) v_1(x, \varepsilon)^\top\right]} \leq \mathbb{E}\left[\norm{v_1(x, \varepsilon) v_1(x, \varepsilon)^\top}\right] \leq \mathbb{E}\left[\norm{v_1(x, \varepsilon)}^2\right].
\end{align*}
Using the bound of $\norm{v_1(x, \varepsilon)}$ in \eqref{equ:inequality-question-01-03:e2-1} and the 4 th moment bound of $\norm{\varepsilon}$, we obtain that
\begin{align}\label{lemma:small-noise-covariance-lower-bound:e4}
    \norm{\Cov{v_1(x, \varepsilon)}} \leq \ell^2 d \mathbb{E}\left[\norm{\varepsilon}^4\right] \leq \ell^2 d \overline{\gamma}.
\end{align}
Note that the norm of $R$ (Equation \eqref{lemma:small-noise-covariance-lower-bound:e3-1}) can be upper bounded by the sum of the norms of the 3 separate terms. Thus, by combining the \eqref{lemma:small-noise-covariance-lower-bound:e3} and \eqref{lemma:small-noise-covariance-lower-bound:e4}, we see that \eqref{lemma:small-noise-covariance-lower-bound:e3-1} holds.
\end{proof}

Lastly, we consider the case when the input of $g$ can be expressed as the sum of a sequence of mutual independent random vectors.

\begin{lemma}\label{lemma:var-sum-of-indep-gaussian}
Let $\{X_i\}_{1\leq i \leq N}$ be a sequence of mean-zero random vectors in $\mathbb{R}^d$ that are mutually independent and satisfies $\underline{\delta} I \preceq \Cov{X_i}$ and $\bE{\norm{X_i}^4} \leq \overline{\gamma}$. Let $g$ be a function that satisfies \eqref{equ:inequality-question-01-03} and \eqref{equ:inequality-question-01-03:e0}. Then, for any positive integer $N$, we have
\begin{align}\label{lemma:var-sum-of-indep-gaussian:e1}
    \Cov{g\left(\sum_{i=1}^N X_i\right)} \succeq N\left(\gamma^2 \underline{\delta} - 2 L \ell d^2 \cdot \overline{\gamma}^{\frac{3}{4}} - \ell^2 d \overline{\gamma}\right) I.
\end{align}
\end{lemma}
\begin{proof}[Proof of Lemma \ref{lemma:var-sum-of-indep-gaussian}]
We use an induction on $N$ to show that \eqref{lemma:var-sum-of-indep-gaussian:e1} holds.

When $N = 1$, \eqref{lemma:var-sum-of-indep-gaussian:e1} holds by setting $x = 0$ and $\varepsilon = X_1$ in Lemma \ref{lemma:small-noise-covariance-lower-bound}.

Suppose \eqref{lemma:var-sum-of-indep-gaussian:e1} holds for $N-1$. Then, for $N$, by the law of total variance, we see that
\begin{align}\label{lemma:var-sum-of-indep-gaussian:e2}
    \Cov{g\left(\sum_{i=1}^N X_i\right)} = \Cov{\mathbb{E}\left[\left.g\left(\sum_{i=1}^N X_i\right)\right| \sum_{i=1}^{N-1} X_i\right]} + \mathbb{E}\left[\Cov{\left.g\left(\sum_{i=1}^N X_i\right)\right| \sum_{i=1}^{N-1} X_i}\right].
\end{align}
For the first term in \eqref{lemma:var-sum-of-indep-gaussian:e2}, we define a new function
\begin{align*}
    \bar{g}(x) \coloneqq \mathbb{E}\left[g(x + X_N)\right].
\end{align*}
Since the random variables $\{X_i\}_{1\leq i\leq N}$ are mutually independent, we observe that
\begin{align*}
    \mathbb{E}\left[\left.g\left(\sum_{i=1}^N X_i\right)\right| \sum_{i=1}^{N-1} X_i\right] = \bar{g}\left[\sum_{i=1}^{N-1} X_i\right].
\end{align*}
One can verify that if $g$ satisfies the conditions in \eqref{equ:inequality-question-01-03} and \eqref{equ:inequality-question-01-03:e0}, then $\bar{g}$ also satisfies the same conditions as $g$ because
\begin{align*}
    \norm{\bar{g}(x) - \bar{g}(x')} = \norm{\mathbb{E}\left[g(x + X_N) - g(x' + X_N)\right]} \leq \mathbb{E}\left[\norm{g(x + X_N) - g(x' + X_N)}\right] \leq L \norm{x - x'}.
\end{align*}
On the other hand, we have
\begin{align*}
    \langle \bar{g}(x) - \bar{g}(x'), x - x'\rangle ={}& \langle \mathbb{E}\left[g(x + X_N) - g(x' + X_N)\right], x - x'\rangle\\
    ={}& \mathbb{E}\left[\langle g(x + X_N) - g(x' + X_N), x - x'\rangle\right] \geq \gamma \norm{x - x'}^2.
\end{align*}
For the Hessian upper/lower bounds, because $\nabla^2 \bar{g}_i(x) = \nabla^2 \mathbb{E}\left[g_i(x + X_N)\right] = \mathbb{E}\left[\nabla^2 g_i(x + X_N)\right]$,
\begin{align*}
    -\ell I \preceq \bar{g}_i(x) \preceq \ell I.
\end{align*}
Therefore, by the induction assumption, we see that
\begin{align}\label{lemma:var-sum-of-indep-gaussian:e3}
    \Cov{\mathbb{E}\left[\left.g\left(\sum_{i=1}^N X_i\right)\right| \sum_{i=1}^{N-1} X_i\right]} ={}& \Cov{\bar{g}\left[\sum_{i=1}^{N-1} X_i\right]}\nonumber\\
    \succeq{}& (N - 1)\left(\gamma^2 \underline{\delta} - 2 L \ell d^2 \cdot \overline{\gamma}^{\frac{3}{4}} - \ell^2 d \overline{\gamma}\right) I.
\end{align}

For the second term in \eqref{lemma:var-sum-of-indep-gaussian:e2}, we note that for any realization $x$ of $\sum_{i=1}^{N-1}X_{i}$, we have 
\begin{align*}
    \Cov{\left.g\left(\sum_{i=1}^N X_i\right)\right| \sum_{i=1}^{N-1} X_i = x} ={}& \Cov{\left.g\left(x + X_N\right)\right| \sum_{i=1}^{N-1} X_i = x} = \Cov{g(x + X_N)}\\
    \succeq{}& \left(\gamma^2 \underline{\delta} - 2 L \ell d^2 \cdot \overline{\gamma}^{\frac{3}{4}} - \ell^2 d \overline{\gamma}\right) I,
\end{align*}
where the conditioning can be removed in the second step because the random variables $\{X_i\}_{1\leq i \leq N}$ are mutually independent, so $g(x + X_N)$ is independent with $\sum_{i=1}^{N-1} X_i$; and we use Lemma \ref{lemma:small-noise-covariance-lower-bound} in the last inequality. Therefore, we obtain that
\begin{align}\label{lemma:var-sum-of-indep-gaussian:e4}
    \mathbb{E}\left[\Cov{\left.g\left(\sum_{i=1}^N X_i\right)\right| \sum_{i=1}^{N-1} X_i}\right] \succeq \left(\gamma^2 \underline{\delta}  - 2 L \ell d^2 \cdot \overline{\gamma}^{\frac{3}{4}} - \ell^2 d \overline{\gamma}\right) I.
\end{align}
Substituting \eqref{lemma:var-sum-of-indep-gaussian:e3} and \eqref{lemma:var-sum-of-indep-gaussian:e4} into \eqref{lemma:var-sum-of-indep-gaussian:e2} shows that \eqref{lemma:var-sum-of-indep-gaussian:e1} still holds for $N$. Thus, we have proved Lemma \ref{lemma:var-sum-of-indep-gaussian} by induction.
\end{proof}

Now we come back to the proof of Theorem \ref{thm:variance-pass-through-function}. By the assumption, we know the distribution of $X$ is identical with the distribution of $\sum_{i=1}^N X_i$, where $X_i$ are i.i.d. random vectors that satisfies $\bE{\norm{X_i}^4} \leq C \cdot N^{-2}$. Thus, we have
\begin{align*}
    \Cov{g(X)} = \Cov{g\left(\sum_{i=1}^N X_i\right)}.
\end{align*}
Note that each $X_i$ satisfies that $\Cov{X_i} = \frac{1}{N}\Cov{X} \succeq \frac{\mu}{N} I$. Applying Lemma \ref{lemma:var-sum-of-indep-gaussian} gives that
\begin{align*}
    \Cov{g(X)} \succeq \left(\mu \gamma^2 - \frac{C^{3/4}}{\sqrt{N}} \cdot 2 L \ell d^2 - \frac{C}{N}\cdot \ell^2 d \overline{\gamma}\right)\cdot I.
\end{align*}
By letting $N$ tends to infinity in the above inequality, we finishes the proof of Theorem \ref{thm:variance-pass-through-function}.

\section{Roadmap to Multi-step Prediction under Well-Conditioned Costs}\label{appendix:multi-step-prediction-well-cond-roadmap}
A limitation of Assumption \ref{assump:one-step-prediction-indep-pairs} in Section \ref{sec:LTV-well-conditioned} is that it only allows the prediction $V_t(\theta)$ to depend on the disturbance $W_t$ at time step $t$. A natural question is whether we can relax the assumption by allowing $V_t(\theta)$ to depend on all future disturbances $W_{t:(T-1)}$. In this section, we present a roadmap towards this generalization and discuss about the potential challenges.

First, we show that the expected cost-to-go function $\mathbb{E}\left[C_t^{\pi^\theta}(x; \Xi)\mid I_t(\theta)\right]$ can be expressed as a function that only depends on the conditional expectations $W_{\tau\mid t}^\theta$ for all $\tau \geq t$, i.e., there exists a function $\tilde{C}_t^{\pi^\theta}$ that satisfies
\begin{align}\label{multi-step-roadmap:e1}
    \tilde{C}_t^{\pi^\theta}(x; W_{t:(T-1)\mid t}^\theta) = \mathbb{E}\left[C_t^{\pi^\theta}(x; \Xi)\mid I_t(\theta)\right].
\end{align}
We show \eqref{multi-step-roadmap:e1} by induction on $t = T, T-1 , \ldots, 0$. Note that the statement holds for $T$. Suppose it holds for $t+1$, by \eqref{equ:auxiliary-cost-to-go-function}, we have
\begin{align*}
    \bar{C}_{t+1}^{\pi^\theta}(x; I_t(\theta)) ={}& \mathbb{E}\left[C_{t+1}^{\pi^\theta}(x + W_t - W_{t\mid t}^\theta; \Xi)\mid I_t(\theta)\right]\\
    ={}& \Condexp{\tilde{C}_{t+1}^{\pi^\theta}(x + W_t - W_{t\mid t}^\theta; W_{(t+1):(T-1)\mid t+1}^\theta)}{I_t(\theta)},
\end{align*}
where we use the induction assumption in the last equation. Define the random variables $\varepsilon_{t\mid t}^\theta \coloneqq W_t - W_{t\mid t}^\theta$ and $\varepsilon_{\tau\mid t}^\theta \coloneqq W_{\tau\mid (t+1)}^\theta - W_{\tau\mid t}^\theta$. Using the properties of joint Gaussian distribution, we know that $\varepsilon_{t:(T-1)\mid t}^\theta$ are independent with $I_t(\theta)$. Therefore,
\begin{align*}
    \bar{C}_{t+1}^{\pi^\theta}(x; I_t(\theta)) ={}& \Condexp{\tilde{C}_{t+1}^{\pi^\theta}(x + \varepsilon_{t\mid t}^\theta; W_{(t+1):(T-1)\mid t}^\theta + \varepsilon_{(t+1):(T-1)\mid t}^\theta)}{I_t(\theta)}\\
    ={}& \mathbb{E}_{\varepsilon_{t:(T-1)\mid t}^\theta}\left[{\tilde{C}_{t+1}^{\pi^\theta}(x + \varepsilon_{t\mid t}^\theta; W_{(t+1):(T-1)\mid t}^\theta + \varepsilon_{(t+1):(T-1)\mid t}^\theta)}\right].
\end{align*}
Thus, $\bar{C}_{t+1}^{\pi^\theta}(x; I_t(\theta))$ can be expressed as a function of $x$ and $W_{(t+1):(T-1)\mid t}^\theta$, and we denote it as
\begin{align}\label{multi-step-roadmap:e2}
    \tilde{\bar{C}}_{t+1}(x; W_{(t+1):(T-1)\mid t}^\theta) \coloneqq \bar{C}_{t+1}^{\pi^\theta}(x; I_t(\theta)).
\end{align}
Therefore, we obtain that
\begin{align*}
    \mathbb{E}\left[C_t^{\pi^\theta}(x; \Xi)\mid I_t(\theta)\right] ={}& h_t^x(x) + (h_t^u\square_{(-B_t)}\bar{C}_{t+1}^{\pi^\theta})(A_t x + W_{t\mid t}^\theta; I_t(\theta))\\
    ={}& h_t^x(x) + (h_t^u\square_{(-B_t)}\tilde{\bar{C}}_{t+1}^{\pi^\theta})(A_t x + W_{t\mid t}^\theta; W_{(t+1):(T-1)\mid t}^\theta).
\end{align*}
Therefore, $\mathbb{E}\left[C_t^{\pi^\theta}(x; \Xi)\mid I_t(\theta)\right]$ can also be expressed in the form $\tilde{C}_t^{\pi^\theta}(x; W_{t:(T-1)\mid t}^\theta)$. Thus, we have shown \eqref{multi-step-roadmap:e1} by induction, with \eqref{multi-step-roadmap:e2} as an intermediate result.

Note that the optimal policy is given by
\begin{align*}
    \pi_t^\theta(x; I_t(\theta)) &\coloneqq \argmin_u \left(h_t^u(u) + \bar{C}_{t+1}^{\pi^\theta}(A_t x + B_t u + W_{t\mid t}^\theta; I_t(\theta))\right)\\
    &= \argmin_u \left(h_t^u(u) + \tilde{\bar{C}}_{t+1}^{\pi^\theta}(A_t x + B_t u + W_{t\mid t}^\theta; W_{(t+1):(T-1)\mid t}^\theta)\right)\\
    &= u_{(h_t^u\square_{-B_t}\tilde{\bar{C}}_{t+1}^{\pi^\theta})}(A_t x + W_{t\mid t}^\theta; W_{(t+1):(T-1)\mid t}^\theta).
\end{align*}
Therefore, by Lemma \ref{lemma:infimal-conv-variant-grad-var-to-solution-var}, we need to establish a covariance lower bound of the gradient
\begin{align*}
\nabla_x\tilde{\bar{C}}_{t+1}^{\pi^\theta}(x + W_{t\mid t}^\theta; W_{(t+1):(T-1)\mid t}^\theta)
\end{align*}
in order to derive a lower bound for the trace of the covariance matrix of $\pi_t^\theta(x; I_t(\theta))$. While this is relatively straightforward when we only have $W_{t\mid t}^\theta$ because it is added directly with $x$, it is much more challenging to also consider the covariance caused by $W_{(t+1):(T-1)\mid t}^\theta$. This is because they affect $\tilde{\bar{C}}_{t+1}^{\pi^\theta}$ through multiple steps of infimal convolutions. Nevertheless, we feel the approach that we describe here is promising if we can derive more properties that are preserved through the infimal convolution operators. We leave this direction as future work.

\end{document}